\documentclass[a4paper,10pt,twoside]{article}
\usepackage{pst-fill,pst-grad}
\usepackage{textcomp}
\usepackage[english]{babel}
\usepackage[utf8x]{inputenc}
\usepackage{graphicx}
\usepackage{amsmath}
\usepackage{float}
\usepackage{fancyhdr}
\usepackage[matrix,arrow,curve]{xy}
\usepackage{pstricks} 
\usepackage{amsmath,amsfonts,verbatim,afterpage,theorem,euscript,mathrsfs,amssymb}
\usepackage{amsfonts}
\usepackage{amssymb}
\usepackage{array}
\usepackage{dsfont}
\usepackage{hyperref}
 
\newtheorem{Definition}{Definition}[section]
\newtheorem{Proposition}{Proposition}[section]
\newtheorem{Lemme}{Lemma}[section]
\newtheorem{Theoreme}{Theorem}
\newtheorem{Corollaire}{Corollary}[section]
\newtheorem{Remarque}{Remark}[section]
\title{\bf Non-local diffusion equations with Lévy-type operators and divergence free drift}
\author{Diego Chamorro\footnote{Laboratoire d'Analyse et de Probabilit\'es, Universit\'e d'Evry Val d'Essonne, 23 Boulevard de France, 91037 Evry Cedex - France, diego.chamorro@univ-evry.fr}}
\begin{document}
\maketitle
\begin{scriptsize}
\abstract{In this paper we are interested in some properties related to the solutions of non-local diffusion equations with divergence free drift. Existence, maximum principle and a positivity principle are proved. In order to study H\"older regularity, we apply a method that relies in the H\"older-Hardy spaces duality and in the molecular characterisation of local Hardy spaces. In these equations, the diffusion is given by Lévy-type operators with an associated Lévy measure satisfying some upper and lower bounds.}\\[3mm]
\textbf{Keywords: Lévy-type operators, Lévy-Khinchin formula, Hölder regularity, molecular Hardy spaces.}
\end{scriptsize}
\section{Introduction}
We study in this article a class of non-local diffusion equations with divergence free drift of the following form:
\begin{equation}\label{Equation0}
\begin{cases}
\partial_t \theta(x,t)-\nabla\cdot(v\,\theta)(x,t)+\mathcal{L}\theta(x,t)=0,\\[5mm]
\theta(x,0)= \theta_0(x),\\[5mm]
\mbox{with }\; div(v)= 0\; \mbox{ and } t\in [0, T].
\end{cases}
\end{equation}
This type of transport-diffusion equations is a generalization of a well-known equation from fluid dynamics. Indeed, in space dimension $n=2$ if $\mathcal{L}=(-\Delta)^{\alpha}$ is the fractional Laplacian, with $0<\alpha\leq 1/2$, and if $v=(-R_2\theta, R_1\theta)$ where $R_{1,2}$ are the Riesz Transforms defined in the Fourier level by $ \widehat{R_j\theta}(\xi)=-\frac{i\xi_j}{|\xi|}\widehat{\theta}(\xi)$ for $j=1,2$, we obtain the quasi-geostrophic equation $(QG)_\alpha$ which has been recently studied by many authors with different approaches and with a variety of results, see \cite{Caffarelli}, \cite{Cordoba}, \cite{KN}, \cite{CW1}, \cite{CW}, \cite{Marchand1} and the references there in for more details.\\ 

Inspired by the work of Kiselev and Nazarov \cite{KN}, it is possible to study the H\"older regularity of the solutions of the $(QG)_{1/2}$ equation by a duality-based method. The aim of this article is to generalize this method to a wider family of operators and we will consider here Lévy-type operators under some hypothesis that will be stated in the lines below. This class of operators corresponds to a natural generalization of recent works where some results are obtained for different operators using quite specific techniques:  for example see the article \cite{MM} where the operator's kernel satisfies some similar bounds to those imposed in our hypothesis.\\

In this paper we will mainly consider problems of existence of the solutions, a maximum principle, a positivity principle and of course we will study H\"older regularity of the solutions of equation (\ref{Equation0}).\\

Let us start by describing our setting in a general way. This framework will be made precise later on.
\begin{itemize}
\item In the formula (\ref{Equation0}) we noted $\theta:\mathbb{R}^n\times [0,T]\longrightarrow \mathbb{R}$ a real-valued function, where $n\geq 2$ is the euclidean dimension. \\
\item The drift (or velocity) term $v$ is such that $v:\mathbb{R}^n\times [0,T]\longrightarrow \mathbb{R}^n$ and we will always assume that $div(v)=0$ and that $v$ belongs to $L^{\infty}([0,T];bmo(\mathbb{R}^n))$. Recall that local $bmo(\mathbb{R}^n)$ space is defined as locally integrable functions $f$ such that
$$\underset{|B|\leq 1}{\sup}\frac{1}{|B|}\int_{B}|f(x)-f_B|dx<M \qquad\mbox{ and }\qquad \underset{|B|>1}{\sup}\frac{1}{|B|}\int_{B}|f(x)|dx<M\qquad \mbox{for a constant } M;$$
we noted $B(R)$ a ball of radius $R>0$ and $f_B=\frac{1}{|B|}\displaystyle{\int_{B(R)}}f(x)dx$. The norm $\|\cdot\|_{bmo}$ is then fixed as the smallest constant $M$ satisfying these two conditions.

\item The operator $\mathcal{L}$ is a Lévy operator which has the following general form called the \textit{Lévy-Khinchin} representation formula:
\begin{equation*}
\mathcal{L}(f)(x)=b\cdot \nabla f(x)+\sum_{j,k=1}^n a_{j,k}\frac{\partial^2 f(x)}{\partial x_j \partial x_k}+\int_{\mathbb{R}^n\setminus \{0\}}\big[f(x)-f(x-y)+y\cdot \nabla f(x)\mathds{1}_{\{|y|\leq 1\}}(y)\big]\Pi(dy),
\end{equation*}
where $b\in \mathbb{R}^n$ is a vector, $a_{j,k}$ are constants (note that the matrix $(a_{j,k})_{1\leq j,k\leq n}$ should be positive semi-definite) and $\Pi$ is a nonnegative Borel measure on $\mathbb{R}^n$ satisfying $\Pi(\{0\})=0$ and
\begin{equation}\label{ConditionGeneraleLevy}
\int_{\mathbb{R}^n}\min(1, |y|^2)\Pi(dy)<+\infty.
\end{equation}
\end{itemize}
In the Fourier level we have $\widehat{\mathcal{L}f\,}(\xi)=a(\xi)\widehat{f}(\xi)$ where the symbol $a(\cdot)$ is given by the Lévy-Khinchin formula
\begin{equation}\label{LevyKhinchine}
a(\xi)=ib\cdot \xi +q(\xi)+\int_{\mathbb{R}^n\setminus \{0\}}\bigg(1-e^{-iy\cdot \xi}-iy\cdot \xi \mathds{1}_{\{|y|<1\}}(y)\bigg)\Pi(dy),\qquad \mbox{where } q(\xi)=\sum_{j,k=1}^na_{j,k}\xi_j\xi_k. 
\end{equation}
Our main references concerning Lévy operators and the Lévy-Khinchin representation formula are the books \cite{Jacob}, \cite{Jacob1} and \cite{Sato}. See also the lecture notes \cite{Karch} for interesting applications to the PDEs.\\

We need to make some assumptions over the Lévy operator considered before. First we will set $b=0$ and $a_{j,k}=0$. We assume then that the measure $\Pi$ is absolutely continuous with respect to the Lebesgue measure, so this measure can be written as $\Pi(dy)=\pi(y)dy$, this hypothesis is important as it simplifies considerably the computations. We will also require some symmetry in the following sense: $\pi(y)=\pi(-y)$. Finally, the most crucial issue concerns estimates over the function $\pi$ and we will assume the inequalities:
\begin{eqnarray}
c_1|y|^{-n-2\alpha}\leq &\pi(y)&\leq c_2|y|^{-n-2\beta} \qquad \mbox{over } |y|\leq 1,\label{DefKernel2}\\[5mm]
0\leq &\pi(y)&\leq c_3|y|^{-n-2\delta} \qquad \mbox{over } |y|> 1,\label{DefKernel3}
\end{eqnarray}
where $c_1,c_2, c_3>0$ are positive constants. We need to define the values of the parameters $\alpha, \beta, \delta$ and we will study the following cases:
\begin{itemize}
\item[\textbf{(a)}] $0<\alpha\leq \beta< 1/2$ and $0<\delta<1/2$,
\item[\textbf{(b)}] $0<\alpha=\beta=\delta< 1/2$,
\item[\textbf{(c)}] $\alpha=\beta=1/2$ and  $0<\delta<1/2$,
\item[\textbf{(d)}] $\alpha=\beta=\delta=1/2$.
\end{itemize}
The choice of these bounds is mainly technical and it will be explained in Remark \ref{RemarqueReg} below.\\

Note that these two conditions (\ref{DefKernel2}) and (\ref{DefKernel3}) imply the next pointwise property which will be useful in the sequel
\begin{equation}\label{EstimationMesure1}
0\leq  \pi(y)\leq c_4(|y|^{-n-2\beta}+|y|^{-n-2\delta})\qquad \mbox{for all }\quad y \in \mathbb{R}^n \mbox{ and } c_4>0. 
\end{equation}

We observe now that these assumptions for the function $\pi$ imply that the operator $\mathcal{L}$ and its symbol $a(\cdot)$ can be rewritten in the following way:
\begin{equation}\label{DefOperator1}
\mathcal{L}(f)(x)=\mbox{v.p.} \int_{\mathbb{R}^n}\big[f(x)-f(x-y)\big]\pi(y)dy
\end{equation}
and 
\begin{equation}\label{LevyKhinchine1}
a(\xi)=\int_{\mathbb{R}^n\setminus \{0\}}\big(1-\cos(\xi\cdot y)\big)\pi(y)dy.
\end{equation}
As we can see, the properties of the operator $\mathcal{L}$ can be easily read, in the real variable or in the Fourier level, by the properties of the function $\pi$. \\

In order to have a better understanding of these properties it is helpful to consider an important example which is given by the fractional Laplacian $(-\Delta)^\alpha$ defined by the expression
$$(-\Delta)^\alpha f(x)=\mbox{v.p.} \int_{\mathbb{R}^n}\frac{f(x)-f(x-y)}{|y|^{n+2\alpha}}dy, \quad \mbox{with }0< \alpha\leq 1/2.$$
Note that we have here $\pi(y)=|y|^{-n-2\alpha}$ and $\pi$ satisfies (\ref{DefKernel2}) and (\ref{DefKernel3}) with $\alpha=\beta=\delta$, so this example corresponds to the cases \textbf{(b)} and \textbf{(d)} stated above. Equivalently, we have a Fourier characterisation by the formula $\widehat{(-\Delta)^\alpha f}(\xi)=|\xi|^{2\alpha}\widehat{f}(\xi)$ so the function $a(\xi)$ is equal to $|\xi|^{2\alpha}$.\\ 

With this example we observe that the lower bound in (\ref{DefKernel2}) guarantees a \textit{diffusion or regularization effect}\footnote{the term ``diffusion'' must be taken in the sense of the PDEs considered by analysts.} like $(-\Delta)^\alpha$ and this is an important assumption for the function $\pi$. Indeed, in some general sense, only the part of the integral (\ref{DefOperator1}) near the origin is critical as $\pi$ satisfies (\ref{DefKernel3}). We note also that the upper bounds given in (\ref{DefKernel2}) and (\ref{DefKernel3}) imply the property (\ref{ConditionGeneraleLevy}) since in any case we have $\beta, \delta \leq 1/2$.
\begin{Remarque}\label{RemarqueReg} As the previous example shows, when $\alpha=\beta=\delta$ we obtain the fractional Laplacian $(-\Delta)^\alpha$ and thus the equation (\ref{Equation0}) studied here can be considered as a linearization of the quasi-geostrophic equation where we have an interesting competition between this operator and the drift term. In the framework of this equation it is classical to distinguish three regimes: \emph{super-critical} if $0<\alpha<1/2$, \emph{critical} if $\alpha=1/2$ and \emph{sub-critical} if $1/2<\alpha<1$, from which only the two first are of interest since in the sub-critical case the regularization effect is in some sense ``stronger'' than the drift, see \cite{CW} for more details.

This explains the upper bound given for the parameters $\alpha,\beta,\delta$. The main reason to divide our study following the cases \textbf{(a)-(d)} is technical as some of the results stated below are valid in some special cases. 
\end{Remarque}
Let us consider more examples: it is shown in Theorem 3.7.7 of \cite{Jacob}, that each continuous negative definite function $a(\cdot)$ can be writen in the form (\ref{LevyKhinchine}), so under hypothesis (\ref{DefKernel2}) and (\ref{DefKernel3}) we can obtain a large class of operators that are in the scope of this work. In the paper \cite{MM} another approach is given: the assumptions for the function $\pi$ are quite similar but they are stated in a different way, furthermore the authors of this article only consider the case $\alpha= \beta=\delta$ in their hypothesis, so our framework is slightly more general. However they allow dependence of the function $\pi$ in the $x$ variable and in the time variable $t$. A further work could follow this path, assuming for example in formula (\ref{DefOperator1}) that $\pi=\pi(x,y,t)$ instead of $\pi=\pi(y)$. Note that some amount of work is already done in this direction, see chapter 4 and Definition 4.5.10 of \cite{Jacob} for more information. 
\subsection*{Presentation of the results}

We assume from now on that the operator $\mathcal{L}$ is of the form (\ref{DefOperator1}). We will work with a function $\pi$ satisfying the hypothesis (\ref{DefKernel2}) and (\ref{DefKernel3}) with the parameters $\alpha,\beta,\delta$ satisfying \textbf{(a)-(d)} unless otherwise specified.\\

In this article we present some results concerning non-local diffusion equation (\ref{Equation0}). Maybe the three first of them are well known for different mathematical communities, so perhaps the only novelty here is the use of the $bmo$ space. Nevertheless we will give the proofs for the sake of completness.

\begin{Theoreme}[Existence and uniqueness for $L^p$ initial data]\label{Theo0}
If $\theta_0\in L^{p}(\mathbb{R}^n)$ with $1\leq p\leq +\infty$ is an initial data, then equation (\ref{Equation0}) has a unique weak solution $\theta\in L^\infty([0,T]; L^p(\mathbb{R}^n))$.
\end{Theoreme}

\begin{Theoreme}[Maximum Principle]\label{Theo1}
Let $\theta_0\in L^{p}(\mathbb{R}^n)$ with $1\leq p\leq +\infty$ be an initial data, then the weak solution of equation (\ref{Equation0}) satisfies the following maximum principle for all $t\in [0, T]$: $\|\theta(\cdot, t)\|_{L^p}\leq \|\theta_0\|_{L^p}$.
\end{Theoreme}

\begin{Theoreme}[Positivity Principle]\label{Theo2}
Let $\beta$ and $\delta$ be the parameters given in cases \textbf{(a)-(d)}. Let $ \frac{n}{2\min(\beta, \delta)}\leq p \leq +\infty$ and $M>0$ a constant, if the initial data $\theta_0\in L^p(\mathbb{R}^n)$ is such that $0\leq \theta_0\leq M$ then the weak solution of equation (\ref{Equation0}) satisfies $0\leq \theta(x,t)\leq M$ for all $t\in [0, T]$. 
\end{Theoreme}
Our main theorem is the following one which is a generalization of a duality method used in the framework of the quasi-geostrophic equation. With this method we obtain a small regularity gain, but for technical reasons we need to consider here the cases \textbf{(c)} and \textbf{(d)}. 
\begin{Theoreme}[Hölder regularity]\label{Theo3}
Let $\mathcal{L}$ be a Lévy operator of the form (\ref{DefOperator1}) with a Lévy measure $\pi$ satisfying hypothesis (\ref{DefKernel2}) and (\ref{DefKernel3}) with $\alpha= \beta=1/2$ and $\delta<1/2$ or $\alpha=\beta=\delta=1/2$. Fix a small time $T_0>0$. Let $\theta_0$ be a function such that $\theta_0\in L^{\infty}(\mathbb{R}^n)$. If $\theta(x,t)$ is a solution for the equation (\ref{Equation0}), then for all time $T_0<t<T$, we have that $\theta(\cdot,t)$ belongs to the H\"older space $\mathcal{C}^{\gamma}(\mathbb{R}^n)$ with $0<\gamma<2\delta<1$ in the case \textbf{(c)} or  $0<\gamma<1$ in the case \textbf{(d)}. 
\end{Theoreme}

The plan of the article is the following: in the section \ref{SecExiUnic} we study existence and uniqueness of solutions with initial data in $L^p$ with $1\leq p<+\infty$.  Section \ref{Sect_PrincipeMax} is devoted to a positivity principle that will be useful in our proofs and section \ref{SecLinfty} studies existence of solution with $\theta_0\in L^\infty$. In section \ref{SeccHolderRegularity} we study the H\"older regularity of the solutions of equation (\ref{Equation0}) by a duality method.
\section{Existence and uniqueness with $L^p$ initial data.}\label{SecExiUnic}
In this section we will study existence and uniqueness for weak solution of equation (\ref{Equation0}) with initial data $\theta_0\in L^p(\mathbb{R}^n)$ where $p\geq 1$. We will start by considering viscosity solutions with an approximation of the velocity field $v$, and we will prove existence and uniqueness for this system. To pass to the limit we will need a further step that is a consequence of the maximum principle. 
\begin{Remarque}
Since the velocity $v$ is a data of the problem, it is equivalent to consider $-v$ instead of $v$, thus for simplicity we fix velocity's sign as in equation (\ref{SistApprox}) below. The same proofs are valid for equation (\ref{Equation0}). 
\end{Remarque}
\subsection{Viscosity solutions}\label{SeccVis}
Before passing to further computations, we give an approximation for functions that belong to the $bmo$ space that will be very useful in the sequel. 
\begin{Lemme}\label{TheoApproxbmo}
Let $f$ be a function in $bmo(\mathbb{R}^n)$. For $k\in \mathbb{N}$, define $f_k$ by
\begin{equation}\label{Forbmoaprox}
f_k(x)=\left\lbrace
\begin{array}{rll}
-k & \mbox{if} & f(x)\leq -k \\[2mm]
f(x) & \mbox{if} & -k\leq f(x)\leq k \\[2mm]
k & \mbox{if} & k\leq f(x).
\end{array}
\right.
\end{equation}
Then $(f_k)_{k\in \mathbb{N}}$ converges weakly to $f$ in $bmo(\mathbb{R}^n)$.
\end{Lemme}
A proof of this lemma can be found in \cite{Stein2}. Having this result in mind, we can begin our study of Theorem \ref{Theo0}. For this, we will work with the following approximation of the equation (\ref{Equation0}):
\begin{equation}\label{SistApprox}
\left\lbrace
\begin{array}{l}
\partial_t \theta(x,t)+ \nabla\cdot(v_{\varepsilon}\;\theta)(x,t)+\mathcal{L}\theta(x,t)=\varepsilon \Delta \theta(x,t)\\[5mm]
\theta(x,0)=\theta_0(x)\\[5mm]
div(v)=0 \quad \mbox{ and } v\in L^{\infty}([0,T];  L^{\infty}(\mathbb{R}^n)).
\end{array}
\right.
\end{equation}
where $v_{\varepsilon}$ is defined by  $v_{\varepsilon}=v\ast \omega_{\varepsilon}$ with $\omega_{\varepsilon}(x)=\varepsilon^{-n}\omega(x/\varepsilon)$ and $\omega\in \mathcal{C}^{\infty}_0(\mathbb{R}^n)$ is a function such that $\displaystyle{\int_{\mathbb{R}^n}}\omega(x)dx=1$. 
Here $\mathcal{L}$ is a Lévy operator of the form (\ref{DefOperator1}) with hypothesis (\ref{DefKernel2}) and (\ref{DefKernel3}) with
$\alpha,\beta,\delta$ satisfying the bounds given in the cases \textbf{(a)-(d)}. Following \cite{Cordoba}, the solutions of this problem are called \textit{viscosity solutions}.\\

Note that the problem (\ref{SistApprox}) admits the following equivalent integral representation:
\begin{equation}\label{FormIntegr}
\theta(x,t)=e^{\varepsilon t\Delta}\theta_0(x)-\int_{0}^t e^{\varepsilon (t-s)\Delta}\nabla \cdot(v_\varepsilon\; \theta)(x,s)ds-\int_{0}^t e^{\varepsilon (t-s)\Delta}\mathcal{L} \theta(x,s)ds,
\end{equation}
In order to prove Theorem \ref{Theo0}, we will first investigate a local result with the following theorem where we will apply the Banach contraction scheme in the space $L^{\infty}([0,T]; L^{p}(\mathbb{R}^n))$ with the norm $\|f\|_{L^\infty (L^p)}=\displaystyle{\underset{t\in [0,T]}{\sup}}\|f(\cdot, t)\|_{L^p}$.
\begin{Theoreme}[Local existence]\label{TheoPointFixe}
Let $1\leq p<+\infty$ and let $\theta_0$ and $v$ be two functions such that $\theta_0\in L^p(\mathbb{R}^n)$, $div(v)=0$ and $v\in L^{\infty}([0,T']; L^{\infty}(\mathbb{R}^n))$. If the initial data satisfies $\|\theta_0\|_{L^p}\leq K$ and if $T'$ is a time small enough, then (\ref{FormIntegr}) has a unique solution $\theta \in L^{\infty}([0,T']; L^{p}(\mathbb{R}^n))$ on the closed ball $\overline{B}(0,2K)\subset L^{\infty}([0,T']; L^{p}(\mathbb{R}^n))$. 
\end{Theoreme}
\begin{Remarque}\label{Rem_Approx}
Observe that we fixed here the velocity $v$ such that $v\in L^{\infty}([0,T'];  L^{\infty}(\mathbb{R}^n))$. This is not very restrictive since by Lemma \ref{TheoApproxbmo} we can construct a sequence $v_k\in L^\infty(\mathbb{R}^n)$ that converge weakly to $v$ in $bmo(\mathbb{R}^n)$.
\end{Remarque}
\textit{\textbf{Proof of Theorem \ref{TheoPointFixe}.}} We note $L_\varepsilon(\theta)$ and  $N^{v}_\varepsilon(\theta)$ the quantities
\begin{eqnarray*}
L_\varepsilon(\theta)(x,t)= \int_{0}^t e^{\varepsilon (t-s)\Delta}\mathcal{L}\theta(x,s)ds & \mbox{ and } & 
N^{v}_\varepsilon(\theta)(x,t)=\int_{0}^t e^{\varepsilon (t-s)\Delta}\nabla \cdot(v_\varepsilon\; \theta)(x,s)ds. \\
\end{eqnarray*}
We begin with general remarks concerning these two formulas. For the first expression we have:
\begin{Proposition}\label{PropoMaj2} If $f\in L^{\infty}([0,T']; L^{p}(\mathbb{R}^n))$, then
\begin{equation}\label{Maj2}
\|L_\varepsilon(f)\|_{L^\infty (L^p)}\leq C \Phi(T',\varepsilon)\; \|f\|_{L^\infty (L^p)}
\end{equation}
where $\Phi(T',\varepsilon)=\left(\frac{T'^{1-\beta}}{\varepsilon^{\beta}}+\frac{T'^{1-\delta}}{\varepsilon^{\delta}}\right)$;  $\left(\frac{T'^{1-\alpha}}{\varepsilon^{\alpha}}\right)$; $\left(\frac{T'^{1/2}}{\varepsilon^{1/2}}+T'+\frac{T'^{1-\delta}}{\varepsilon^{\delta}}\right)$ and $\left(\frac{T'^{1/2}}{\varepsilon^{1/2}}\right)$, for the cases \textbf{(a)-(d)} respectively.
\end{Proposition}
\textit{\textbf{Proof.}}
We write
$$\|L_\varepsilon(f)\|_{L^\infty (L^p)}=\underset{0<t<T'}{\sup} \left\|\int_{0}^t e^{\varepsilon (t-s)\Delta}\mathcal{L} f(\cdot,s)ds\right\|_{L^p}=\underset{0<t<T'}{\sup} \left\|\int_{0}^t \mathcal{L} f\ast h_{\varepsilon (t-s)}(\cdot,s)ds\right\|_{L^p}$$
where $h_t$ is the heat kernel on $\mathbb{R}^n$. By the properties of the Lévy operator $\mathcal{L}$ we can write $\mathcal{L}f\ast h_{\varepsilon(t-s)}=f\ast \mathcal{L}h_{\varepsilon(t-s)}$ and then we obtain the estimate
\begin{equation}\label{Phi1}
\|L_\varepsilon(f)\|_{L^\infty (L^p)}\leq \underset{0<t<T'}{\sup} \displaystyle{\int_{0}^t} \|f(\cdot, s)\|_{L^p} \|\mathcal{L}h_{\varepsilon(t-s)}\|_{L^1}ds\leq  \|f\|_{L^\infty (L^p)}\underset{0<t<T'}{\sup} \int_{0}^t  \|\mathcal{L}h_{\varepsilon(t-s)}\|_{L^1}ds.
\end{equation}
We need now to study the quantity $\|\mathcal{L}h_{\varepsilon(t-s)}\|_{L^1}$, for this we will use Besov spaces and a short lemma. We recall that for $0<s<1$ and $1\leq p <+\infty$, homogeneous Besov spaces $\dot{B}^{s,p}_p(\mathbb{R}^n)$ may be defined as
$$\|f\|_{\dot{B}^{s,p}_p}=\left(\int_{\mathbb{R}^n} \int_{\mathbb{R}^n} \frac{|f(x)-f(x-y)|^p}{|y|^{n+ps}}dydx\right)^{1/p}.$$
Now, here is the lemma:
\begin{Lemme}\label{LemmeEstimateOperK}
Let $\mathcal{L}$ be a Lévy operator satisfying the hypothesis stated above.
\begin{itemize}
\item[\textbf{(a)}] If $0<\alpha\leq\beta<1/2$ and $0<\delta<1/2$ then, for all $f\in \dot{B}^{2\beta,1}_1(\mathbb{R}^n)\cap \dot{B}^{2\delta,1}_1(\mathbb{R}^n)$ we have $\|\mathcal{L}f\|_{L^1}\leq \|f\|_{\dot{B}^{2\beta,1}_1}+\|f\|_{\dot{B}^{2\delta,1}_1}$. 
In particular we have for the heat kernel $\|\mathcal{L}h_{\varepsilon(t-s)}\|_{L^1}\leq C\big([\varepsilon(t-s)]^{-\beta}+[\varepsilon(t-s)]^{-\delta}\big)$.\\
\item[\textbf{(b)}] If $\alpha=\beta=\delta < 1/2$, we have $\mathcal{L}=(-\Delta)^{\alpha}$ and thus $\|\mathcal{L}h_{\varepsilon(t-s)}\|_{L^1}\leq C[\varepsilon(t-s)]^{-\alpha}$.\\
\item[\textbf{(c)}] If $\alpha=\beta=1/2$ and $0<\delta<1/2$ we have $\|\mathcal{L}f\|_{L^1}\leq C\big(\|(-\Delta)^{1/2}f\|_{L^1}+\|f\|_{L^1}+\|f\|_{\dot{B}^{2\delta,1}_1}\big)$ where the quantities above are assumed to be bounded. In particular we have $\|\mathcal{L}h_{\varepsilon(t-s)}\|_{L^1}\leq C\big([\varepsilon(t-s)]^{-1/2}+1+[\varepsilon(t-s)]^{-\delta}\big)$.\\
\item[\textbf{(d)}] If $\alpha=\beta=\delta=1/2$, we have $\mathcal{L}=(-\Delta)^{1/2}$ and thus $\|\mathcal{L}h_{\varepsilon(t-s)}\|_{L^1}\leq C[\varepsilon(t-s)]^{-1/2}$.
\end{itemize}
\end{Lemme}
\textbf{\textit{Proof of the lemma.}} By homogeneity the cases \textbf{(b)} and \textbf{(d)} are straightforward. 
If $0<\alpha\leq\beta<1/2$ and $0<\delta<1/2$, using (\ref{DefKernel2}) and (\ref{DefKernel3}) we obtain
\begin{eqnarray*}
\|\mathcal{L}f\|_{L^1}&\leq &\int_{\mathbb{R}^n} \int_{\mathbb{R}^n} \frac{|f(x)-f(x-y)|}{|y|^{n+2\beta}}dydx+ \int_{\mathbb{R}^n} \int_{\mathbb{R}^n} \frac{|f(x)-f(x-y)|}{|y|^{n+2\delta}}dydx= \|f\|_{\dot{B}^{2\beta,1}_1}+\|f\|_{\dot{B}^{2\delta,1}_1}.
\end{eqnarray*}
If $\alpha=\beta=1/2$ and $\delta<1/2$, we simply write 
\begin{eqnarray*}
\|\mathcal{L}f\|_{L^1}&\leq&\int_{\mathbb{R}^n} \left|\mbox{v.p.}\int_{\{|y|\leq 1\}}\big[f(x)-f(x-y)\big]\pi(y)dy\right|dx+\int_{\mathbb{R}^n}\left|\int_{\{|y|> 1\}}[f(x)-f(x-y)]\pi(y)dy\right|dx\\
&\leq&\int_{\mathbb{R}^n} \left|\mbox{v.p.}\int_{\{|y|\leq 1\}}\frac{f(x)-f(x-y)}{|y|^{n+1}}dy\right|dx+\|f\|_{\dot{B}^{2\delta,1}_1}.
\end{eqnarray*}
Now, since $(-\Delta)^{1/2}f(x)=\mbox{v.p.}\displaystyle{\int_{\mathbb{R}^n}}\frac{f(x)-f(x-y)}{|y|^{n+1}}dy$ it is easy to obtain that 
\begin{equation*}
\int_{\mathbb{R}^n} \left|\mbox{v.p.}\int_{\{|y|\leq 1\}}\frac{f(x)-f(x-y)}{|y|^{n+1}}dy\right|dx\leq \|(-\Delta)^{1/2}f\|_{L^1}+C\|f\|_{L^1}.
\end{equation*}
Finally, by homogeneity and since the heat kernel $h_t$ is a smooth function, we obtain the wished estimates for this function. \hfill$\blacksquare$\\

With these estimates at our disposal for the quantity $ \|\mathcal{L}h_{\varepsilon(t-s)}\|_{L^1}$, we obtain for (\ref{Phi1}) -after an integration in time and following the different cases- the inequality $\|L_\varepsilon(f)\|_{L^\infty (L^p)}\leq C \Phi(T',\varepsilon)\|f\|_{L^\infty (L^p)}$ and the Proposition \ref{PropoMaj2} is proven. \hfill$\blacksquare$\\

For the term $N^{v}_\varepsilon$ we have:
\begin{Proposition} If $f\in L^{\infty}([0,T']; L^{p}(\mathbb{R}^n))$ and if $v\in  L^{\infty}([0,T']; L^{\infty}(\mathbb{R}^n)) $, then 
\begin{equation}\label{Maj3}
\|N^v_\varepsilon(f)\|_{L^\infty (L^p)}\leq C \sqrt{\frac{T'}{\varepsilon}}\;\|v\|_{L^\infty (L^\infty)} \|f\|_{L^\infty (L^p)}
\end{equation}
\end{Proposition}
\textit{\textbf{Proof.}} We write:
\begin{eqnarray*}
\|N^v_\varepsilon(f)\|_{L^\infty (L^p)}& =&\underset{0<t<T'}{\sup} \left\|\int_{0}^t e^{\varepsilon (t-s)\Delta} \nabla \cdot(v_\varepsilon f)(\cdot,s)ds\right\|_{L^p}=\underset{0<t<T'}{\sup} \left\|\int_{0}^t \nabla \cdot(v_\varepsilon f)\ast h_{\varepsilon (t-s)}(\cdot,s)ds\right\|_{L^p}\\
&\leq & \underset{0<t<T'}{\sup}\int_{0}^t  \left\|v_\varepsilon f(\cdot,s)\right\|_{L^p} \left\|\nabla h_{\varepsilon (t-s)}\right\|_{L^1} ds
\leq  \underset{0<t<T'}{\sup}\int_{0}^t  \left\|v_\varepsilon(\cdot,s)\right\|_{L^\infty}  \left\|f(\cdot,s)\right\|_{L^p} C(\varepsilon(t-s))^{-1/2} ds\\
&\leq &  \|v\|_{L^\infty (L^\infty)}\left\|f\right\|_{L^\infty (L^p)}  \underset{0<t<T'}{\sup}\int_{0}^t C(\varepsilon(t-s))^{-1/2} ds\leq C \sqrt{\frac{T'}{\varepsilon}}  \|v\|_{L^\infty (L^\infty)} \left\|f\right\|_{L^\infty (L^p)}.
\end{eqnarray*}
\hfill$\blacksquare$\\

To finish the preliminary remarks we note, that since $e^{\varepsilon t \Delta}$ is a contraction operator, the estimate $\|e^{\varepsilon t \Delta}f\|_{L^p}\leq \|f\|_{L^p}$ is valid for all function $f\in L^{p}(\mathbb{R}^n)$ with $1\leq p\leq +\infty$, for all $t>0$ and all $\varepsilon>0$. Thus, we have
\begin{equation}\label{Maj1}
\|e^{\varepsilon t \Delta}f\|_{L^\infty (L^p)}\leq \|f\|_{L^p}.
\end{equation}
Now we can use the Banach contraction scheme: we construct a sequence of functions in the following way
$$\theta_{n+1}(x,t)=e^{\varepsilon t\Delta}\theta_0(x)-L_{\varepsilon}(\theta_n)(x,t)-N^v_{\varepsilon}(\theta_n)(x,t)$$
and we take the $L^\infty L^p$-norm of this expression to obtain
$$ \|\theta_{n+1}\|_{L^\infty (L^p)}\leq\|e^{\varepsilon t\Delta}\theta_0\|_{L^\infty (L^p)}+\|L_{\varepsilon}(\theta_n)\|_{L^\infty (L^p)}+\|N^v_{\varepsilon}(\theta_n)\|_{L^\infty (L^p)}$$
Using estimates (\ref{Maj2}), (\ref{Maj3}) and (\ref{Maj1}) we have
$$ \|\theta_{n+1}\|_{L^\infty (L^p)}\leq \|\theta_0\|_{L^p}+C\left( \Phi(T', \varepsilon)+\frac{T'^{1/2}}{\varepsilon^{1/2}}\|v\|_{L^\infty (L^\infty)} \right)\|\theta_n\|_{L^\infty (L^p)}$$
Thus, if $\|\theta_0\|_{L^p}\leq K$ and if we define the time $T'$ to be such that $C\left( \Phi(T', \varepsilon)+\frac{T'^{1/2}}{\varepsilon^{1/2}}\|v\|_{L^\infty (L^\infty)} \right)\leq 1/2$, 
we have by iteration that $\|\theta_{n+1}\|_{L^\infty (L^p)}\leq 2 K$: the sequence $(\theta_n)_{n\in \mathbb{N}}$ constructed from initial data $\theta_0$ belongs to the closed ball $\overline{B}(0, 2K)$. In order to finish this proof, let us show that $\theta_n \longrightarrow \theta$ in $L^{\infty}([0,T']; L^{p}(\mathbb{R}^n))$. For this we write
$$\|\theta_{n+1}-\theta_n\|_{L^\infty (L^p)}\leq \|L_\varepsilon(\theta_{n}-\theta_{n-1})\|_{L^\infty (L^p)}+\|N^v_\varepsilon(\theta_{n}-\theta_{n-1})\|_{L^\infty (L^p)}$$
and using the previous results we have
$$\|\theta_{n+1}-\theta_n\|_{L^\infty (L^p)}\leq C\left(\Phi(T', \varepsilon)+\frac{T'^{1/2}}{\varepsilon^{1/2}}\|v\|_{L^\infty (L^\infty)} \right)\|\theta_{n}-\theta_{n-1}\|_{L^\infty (L^p)}$$
so, by iteration we obtain
$$\|\theta_{n+1}-\theta_n\|_{L^\infty (L^p)}\leq \left[C\left(\Phi(T', \varepsilon)+\frac{T'^{1/2}}{\varepsilon^{1/2}}\|v\|_{L^\infty (L^\infty)} \right)\right]^{n}\|\theta_{1}-\theta_0\|_{L^\infty (L^p)}$$
hence, with the definition of $T'$ it comes
$\|\theta_{n+1}-\theta_n\|_{L^\infty (L^p)}\leq \left(\frac{1}{2}\right)^{n}\|\theta_{1}-\theta_0\|_{L^\infty (L^p)}$.
Finally, if $n\longrightarrow +\infty$, the sequence $(\theta_n)_{n\in \mathbb{N}}$ convergences towards $\theta$ in $L^\infty([0,T'];L^p(\mathbb{R}^n))$. Since it is a Banach space we deduce uniqueness for the solution $\theta$ of problem (\ref{FormIntegr}). The proof of Theorem \ref{TheoPointFixe} is finished.\hfill$\blacksquare$\\

\begin{Corollaire}\label{CorDepContinue}
The solution constructed above depends continuously on the initial value $\theta_0$.
\end{Corollaire}
\textit{\textbf{Proof.}}
Let $\varphi_0, \theta_0\in L^{p}(\mathbb{R}^n)$ be two initial values and let $\varphi$ and $\theta$ be the associated solutions. We write
$$\theta(x,t)-\varphi(x,t)=e^{\varepsilon t\Delta}(\theta_0(x)-\varphi_0(x))-L_\varepsilon(\theta-\varphi)(x,t)-N^v_\varepsilon(\theta-\varphi)(x,t)$$
Taking $L^\infty L^p$-norm in formula above and applying the same previous calculations one obtains
\begin{equation*}
\|\theta-\varphi\|_{L^\infty (L^p)}\leq  \|\theta_0-\varphi_0\|_{L^p}+ C_0\|\theta-\varphi\|_{L^\infty (L^p)}
\end{equation*}
This shows continuous dependence of the solution since $C_0=C\left(\Phi(T', \varepsilon)+\frac{T'^{1/2}}{\varepsilon^{1/2}}\|v\|_{L^\infty (L^\infty)} \right)\leq 1/2$.\hfill$\blacksquare$
\begin{Remarque}[From Local to Global]
Once we obtain a local result, global existence easily follows by a simple iteration since problems studied here (equations (\ref{Equation0}) or (\ref{SistApprox})) are linear as the velocity $v$ does not depend on $\theta$.
\end{Remarque}
We study now the regularity of the solutions constructed by this method.
\begin{Theoreme} Solutions of the approximated problem (\ref{SistApprox}) are smooth.
\end{Theoreme}
\textit{\textbf{Proof}.}
By iteration we will prove that $\theta \in \displaystyle{\bigcap_{0<T_0<T_1<t<T_2<T^\ast}}L^\infty([0,t]; \dot{W}^{\frac{k}{2},p}(\mathbb{R}^n))$ for all $k\geq 0$ where we define the Sobolev space $\dot{W}^{s,p}(\mathbb{R}^n)$ for $ s\in \mathbb{R}$ and $1<p<+\infty$ by  $\|f\|_{\dot{W}^{s,p}}=\|(-\Delta)^{s/2}f\|_{L^p}$. Note that this is true for $k=0$. So let us assume that it is also true for $k>0$ and we will show that it is still true for $k+1$.\\

Set $t$ such that $0<T_0<T_1<t<T_2<T^\ast$ and let us consider the next problem
$$\theta(x,t)=e^{\varepsilon (t-T_0)\Delta}\theta(x, T_0)-\int_{T_0}^t e^{\varepsilon (t-s)\Delta}\nabla \cdot(v_\varepsilon\; \theta)(x,s)ds-\int_{T_0}^t e^{\varepsilon (t-s)\Delta}\mathcal{L} \theta(x,s)ds
$$
We have then the following estimate
\begin{eqnarray*}
\|\theta\|_{L^\infty (\dot{W}^{\frac{k+1}{2},p})}&\leq& \|e^{\varepsilon (t-T_0)\Delta}\theta(\cdot, T_0)\|_{L^\infty (\dot{W}^{\frac{k+1}{2},p})}\\[5mm]
& &+\left\|\int_{T_0}^t e^{\varepsilon (t-s)\Delta}\nabla \cdot(v_\varepsilon\; \theta)(\cdot,s)ds\right\|_{L^\infty (\dot{W}^{\frac{k+1}{2},p})}
+\left\|\int_{T_0}^t e^{\varepsilon (t-s)\Delta}\mathcal{L}\theta(\cdot,s)ds\right\|_{L^\infty (\dot{W}^{\frac{k+1}{2},p})}
\end{eqnarray*}
Now, we will treat separately each of the previous terms. 
\begin{enumerate}
\item[(i)] For the first one we have
$$\|e^{\varepsilon (t-T_0)\Delta}\theta(\cdot, T_0)\|_{\dot{W}^{\frac{k+1}{2},p}}=\|\theta(\cdot, T_0)\ast(-\Delta)^{\frac{k+1}{4}}h_{\varepsilon (t-T_0)}\|_{L^p}\leq\|\theta(\cdot, T_0)\|_{L^p}\| (-\Delta)^{\frac{k+1}{4}}h_{\varepsilon (t-T_0)}\|_{L^1} $$ 
where $h_t$ is the heat kernel, so we can write
\begin{equation*}
\|e^{\varepsilon (t-T_0)\Delta}\theta(\cdot, T_0)\|_{L^\infty (\dot{W}^{\frac{k+1}{2},p})}\leq C\|\theta(\cdot, T_0)\|_{L^p}\sup\left\{ \left[\varepsilon (t-T_0)\right]^{- \frac{k+1}{4}}; 1\right\}
\end{equation*}
\item[(ii)] For the second term, one has
\begin{eqnarray*}
\left\|\int_{T_0}^t e^{\varepsilon (t-s)\Delta}\nabla \cdot(v_\varepsilon\; \theta)(\cdot,s)ds\right\|_{\dot{W}^{\frac{k+1}{2},p}}&\leq &\int_{T_0}^t\|\nabla  \cdot (v_{\varepsilon}\;\theta)\ast h_{\varepsilon(t-s)}\|_{\dot{W}^{\frac{k+1}{2},p}}ds\\
&\leq &\int_{T_0}^t\|(-\Delta)^{\frac{k+1}{4}}\big[\nabla  \cdot (v_{\varepsilon}\;\theta)\ast h_{\varepsilon(t-s)}\big]\|_{L^p}ds\\
&\leq &C\int_{T_0}^t\|v_\varepsilon\;\theta(\cdot,s)\|_{\dot{W}^{\frac{k}{2},p}} \left[\varepsilon (t-s)\right]^{- \frac{3}{4}} ds.
\end{eqnarray*}
Note now that we have here the estimations below for $N\geq k/2$
\begin{eqnarray*}
\|v_\varepsilon\theta(\cdot,s)\|_{\dot{W}^{\frac{k}{2},p}}&\leq& \|v_\varepsilon(\cdot,s)\|_{\mathcal{C}^N} \|\theta(\cdot,s)\|_{\dot{W}^{\frac{k}{2},p}}\leq C \varepsilon^{-N}\|v(\cdot,s)\|_{L^\infty}\|\theta(\cdot,s)\|_{\dot{W}^{\frac{k}{2},p}}
\end{eqnarray*}
hence, we can write
\begin{eqnarray*}
\left\|\int_{T_0}^t e^{\varepsilon (t-s)\Delta}\nabla \cdot(v_\varepsilon\; \theta)(\cdot,s)ds\right\|_{L^\infty (\dot{W}^{\frac{k+1}{2},p})}&\leq &C \|v\|_{L^\infty (L^\infty)} \|\theta\|_{L^\infty (\dot{W}^{\frac{k}{2},p})}\int_{T_0}^t  \varepsilon^{-N}\sup\left\{  \left[\varepsilon (t-s)\right]^{- \frac{3}{4}} ;1\right\}ds
\end{eqnarray*}
\item[(iii)] Finally, for the last term we have
\begin{eqnarray*}
\left\|\int_{T_0}^t e^{\varepsilon (t-s)\Delta}\mathcal{L}\theta(\cdot,s)ds\right\|_{\dot{W}^{\frac{k+1}{2},p}} &\leq & \int_{T_0}^t \left\|(-\Delta)^{\frac{k}{4}}\theta(\cdot,s)\ast \mathcal{L}(-\Delta)^{\frac{1}{4}}h_{\varepsilon(t-s)}\right\|_{L^{p}}ds\\
&\leq &\int_{T_0}^t \left\|\theta(\cdot,s)\right\|_{\dot{W}^{\frac{k}{2},p}}\|\mathcal{L}(-\Delta)^{\frac{1}{4}}h_{\varepsilon(t-s)}\|_{L^{1}}ds\end{eqnarray*}
now, applying Lemma \ref{LemmeEstimateOperK} to the function $(-\Delta)^{\frac{1}{4}}h_{\varepsilon(t-s)}$ we obtain by homogeneity that 
$$\|\mathcal{L}(-\Delta)^{\frac{1}{4}}h_{\varepsilon(t-s)}\|_{L^{1}}\leq \phi(\varepsilon(t-s))$$ 
 where $\phi(\varepsilon(t-s))=\big([\varepsilon (t-s)]^{-\frac{1+4\beta}{4}}+[\varepsilon (t-s)]^{-\frac{1+4\delta}{4}} \big)$; $\big([\varepsilon (t-s)]^{-\frac{1+4\alpha}{4}}\big)$; $\big([\varepsilon (t-s)]^{-\frac{3}{4}}+[\varepsilon (t-s)]^{-\frac{1}{4}}+[\varepsilon(t-s)]^{-\frac{1+4\delta}{4}} \big)$ and $\big([\varepsilon (t-s)]^{-\frac{3}{4}}\big)$ for the cases \textbf{(a)-(d)} respectively. So we obtain:
\begin{eqnarray*}
\left\|\int_{T_0}^t e^{\varepsilon (t-s)\Delta} \mathcal{L}  \theta(\cdot,s)ds\right\|_{L^\infty (\dot{W}^{\frac{k+1}{2},p})}&\leq & C\|\theta \|_{L^{\infty}(\dot{W}^{\frac{k}{2},p})}\int_{T_0}^t \sup\left\{\phi(\varepsilon(t-s));1\right\}ds.
\end{eqnarray*}
\end{enumerate}
Now, with formulas (i)-(iii) at our disposal, we have that the norm $\|\theta\|_{L^\infty (\dot{W}^{\frac{k+1}{2},p})}$ is controlled for all $\varepsilon>0$: we have proven spatial regularity.\\

Time regularity follows since we have
$$\frac{ \partial^k}{\partial t^k}\theta(x,t)+\nabla \cdot \left(\frac{ \partial^k}{\partial t^k} (v_\varepsilon\,\theta)\right)(x,t)+\mathcal{L}\left(\frac{ \partial^k}{\partial t^k} \theta\right)(x,t)=\varepsilon \Delta \left(\frac{ \partial^k}{\partial t^k} \theta\right)(x,t).$$
\begin{flushright}$\blacksquare$\end{flushright}
\begin{Remarque}\label{RemarkEpsilonDep}
The solutions $\theta(\cdot,\cdot)$ constructed above depend on $\varepsilon$. 
\end{Remarque}
\subsection{Maximum principle and Besov regularity} 

The maximum principle we are studying here will be a consequence of few inequalities, some of them are well known. We will start with the solutions $\theta(\cdot, \cdot)$ obtained in the previous section:
\begin{Proposition}[Viscosity version of Theorem \ref{Theo1}]\label{PropoViscosityMaxPrinc}
Let $\theta_0\in L^{p}(\mathbb{R}^n)$ with $1\leq p\leq +\infty$ be an initial data, then the associated solution of the viscosity problem  (\ref{SistApprox}) satisfies the following maximum principle for all $t\in [0, T]$: $\|\theta(\cdot, t)\|_{L^p}\leq \|\theta_0\|_{L^p}$.
\end{Proposition}
\textit{\textbf{Proof.}} We write for $1\leq p<+\infty$:
\begin{eqnarray*}
\frac{d}{dt}\|\theta(\cdot, t)\|^p_{L^p}&=&p\int_{\mathbb{R}^n}|\theta|^{p-2}\theta\bigg(\varepsilon \Delta \theta-\nabla \cdot (v_\varepsilon \,\theta)-\mathcal{L}\theta\bigg)dx=p\varepsilon\int_{\mathbb{R}^n}|\theta|^{p-2}\theta\Delta \theta dx-p\int_{\mathbb{R}^n}|\theta|^{p-1}sgn(\theta) \mathcal{L}\theta dx\\
\end{eqnarray*}
where we used the fact that $div(v)=0$. Thus, we have
$$\frac{d}{dt}\|\theta(\cdot, t)\|^p_{L^p}-p\varepsilon\int_{\mathbb{R}^n}|\theta|^{p-2}\theta \Delta \theta dx+p\int_{\mathbb{R}^n}|\theta|^{p-1}sgn(\theta)\mathcal{L}\theta dx=0,$$
and integrating in time we obtain
\begin{equation}\label{Form1}
\|\theta(\cdot, t)\|^p_{L^p}-p\varepsilon\int_{0}^t\int_{\mathbb{R}^n}|\theta|^{p-2}\theta\Delta \theta dxds+p\int_0^t\int_{\mathbb{R}^n}|\theta|^{p-1}sgn(\theta)\mathcal{L} \theta dxds=\|\theta_0\|^p_{L^p}.
\end{equation}
To finish, we have the following lemma
\begin{Lemme} 
The quantities $-p\varepsilon\displaystyle{\int_{\mathbb{R}^n}}|\theta|^{p-2}\theta \Delta \theta dx$ and $p\displaystyle{\int_0^t\int_{\mathbb{R}^n}}|\theta|^{p-1}sgn(\theta) \mathcal{L} \theta dxds$ are both positive.
\end{Lemme}
\textit{\textbf{Proof.}} For the first expression, since $e^{\varepsilon s\Delta}$  is a contraction semi-group we have
 $\|e^{\varepsilon s\Delta}f\|_{L^p}\leq \|f\|_{L^p}$ for all $s>0$ and all $f\in L^p(\mathbb{R}^n)$. Thus $F(s)=\|e^{\varepsilon s\Delta}f\|_{L^p}$ is decreasing in $s$; taking the derivative in $s$ and evaluating in $s=0$ we obtain the desired result. The positivity of the second expression follows immediately from the \textit{Strook-Varopoulos estimate} for general Lévy-type operators given by the following formula (see remark 1.23 of \cite{Karch} for a proof, more details can be found in \cite{Strook} and \cite{Varopoulos}):
\begin{equation} \label{Strook-Varopoulos}
C\langle \mathcal{L}|\theta|^{p/2},|\theta|^{p/2}  \rangle\leq \langle \mathcal{L}\theta, |\theta|^{p-1}sgn(\theta) \rangle
\end{equation} 
To conclude it is enough to note that $\langle \mathcal{L}|\theta|^{p/2},|\theta|^{p/2}\rangle=\|\mathcal{L}^{\frac{1}{2}}|\theta|^{p/2}\|_{L^2}^2 \geq 0$, where the operator $\mathcal{L}^{\frac{1}{2}}$ is defined by the formula $(\mathcal{L}^{\frac{1}{2}}f)^{\widehat{\quad}}(\xi)=a^{\frac{1}{2}}(\xi)\widehat{f}(\xi)$. \hfill$\blacksquare$\\

Getting back to (\ref{Form1}), we have that all these quantities are bounded and positive and we write for all $1\leq p<+\infty$:
$$\|\theta(\cdot, t)\|_{L^p}\leq \|\theta_0\|_{L^p}.$$
Since $\|\theta(\cdot, t)\|_{L^p}\underset{p\to +\infty}{\longrightarrow}\|\theta(\cdot, t)\|_{L^\infty}$, the maximum principle is proven for viscosity solutions.  \hfill$\blacksquare$\\

In order to deal with Theorem \ref{Theo1} we will need some aditional results. Indeed, a more detailed study of expression (\ref{Form1}) above will lead us to a result concerning weak solution's regularity. 
\begin{Lemme}\label{LemmaSymbol}
If the function $\pi$ satisfies the conditions (\ref{DefKernel2}) and (\ref{DefKernel3}), then we have for the cases \textbf{(a)-(d)} the following pointwise estimates on the symbol $a(\cdot)$ for all $\xi \in \mathbb{R}^n$:
\begin{enumerate}
\item[1)] $a(\xi)\leq |\xi|^{2\beta}+|\xi|^{2\delta}$
\item[2)] $|\xi|^{2\alpha}\leq a(\xi)+C$. 
\end{enumerate}
\end{Lemme}
\textit{\textbf{Proof.}} We use the Lévy-Khinchin formula to obtain  $|\xi|^{2\alpha}=\displaystyle{\int_{\mathbb{R}^n\setminus\{0\}}}\big(1-\cos(y\cdot \xi)\big)|y|^{-n-2\alpha}dy$. It is enough to apply the hypothesis (\ref{DefKernel2}), (\ref{DefKernel3}) and to use the inequality (\ref{EstimationMesure1}) to conclude.\hfill $\blacksquare$\\

\begin{Theoreme}[Besov Regularity]\label{TheoBesov} Let $\mathcal{L}$ be a Lévy-type operator of the form (\ref{DefOperator1}) with hypothesis (\ref{DefKernel2}) and (\ref{DefKernel3}) for the measure $\pi$ with $\alpha,\beta,\delta$ satisfying the bounds given in the cases \textbf{(a)-(d)}. Let $2\leq p <+\infty$ and let $f:\mathbb{R}^n\longrightarrow \mathbb{R}$ be a function such that $f\in L^{p}(\mathbb {R}^n)$ and
$$\int_{\mathbb{R}^n}|f(x)|^{p-2}f(x)\mathcal{L} f(x)dx<+\infty, \quad\mbox{ then  }\quad f\in \dot{B}^{2\alpha/p,p}_{p}(\mathbb{R}^n).$$
\end{Theoreme}
\textit{\textbf{Proof}.} We will prove the following estimates valid for a positive function $f$:
\begin{equation}\label{FormuleEstima1}
\|f\|^p_{ \dot{B}^{2\alpha/p,p}_{p}}\leq C\|f^{p/2}\|^2_{\dot{B}^{\alpha,2}_{2}}\leq \|f^{p/2}\|^2_{L^2} +\int_{\mathbb{R}^n}|f(x)|^{p-2}f(x)\mathcal{L} f(x)dx
\end{equation}
The first inequality can be found in \cite{PGDCH}, so we only need to focus on the right-hand side of the previous estimate. For this, we will start assuming that the function $f$ is positive. \\

Using Plancherel's formula, the characterisation of $\mathcal{L}^{\frac{1}{2}}$ via the symbol $a^{\frac{1}{2}}(\xi)$ and Lemma \ref{LemmaSymbol}  we write
\begin{eqnarray*}
\|f^{p/2}\|^2_{\dot{B}^{\alpha,2}_{2}}&=&\|f^{p/2}\|_{\dot{H}^\alpha}^2=\int_{\mathbb{R}^n}|\xi|^{2\alpha}|\widehat{f^{p/2}}(\xi)|^2d\xi\leq  \int_{\mathbb{R}^n}(a^{\frac{1}{2}}(\xi)+C)^2|\widehat{f^{p/2}}(\xi)|^2d\xi \leq  c\left(\|f^{p/2}\|^2_{L^2}+ \|\mathcal{L}^{\frac{1}{2}}f^{p/2}\|^2_{L^2}\right).
\end{eqnarray*}
Now, using the  Strook-Varopoulos inequality (\ref{Strook-Varopoulos}) we have
$$\|f^{p/2}\|^2_{L^2}+\|\mathcal{L}^{\frac{1}{2}}f^{p/2}\|^2_{L^2}\leq  \|f^{p/2}\|^2_{L^2}+ c\int_{\mathbb{R}^n} f^{p-2}f \mathcal{L}f dx$$
So inequality (\ref{FormuleEstima1}) is proven for positive functions. For the general case we write $f(x)=f_+(x)-f_-(x)$ where $f_\pm(x)$ are positive functions with disjoint support and we have:
\begin{eqnarray}\label{FormDecomp}
\int_{\mathbb{R}^n}|f(x)|^{p-2}f(x)\mathcal{L} f(x)dx&=&\int_{\mathbb{R}^n}f_+(x)^{p-2}f_+(x)\mathcal{L} f_+(x)dx+\int_{\mathbb{R}^n}f_-(x)^{p-2}f_-(x)\mathcal{L} f_-(x)dx\\
&-&\int_{\mathbb{R}^n}f_+(x)^{p-2}f_+(x)\mathcal{L} f_-(x)dx-\int_{\mathbb{R}^n}f_-(x)^{p-2}f_-(x)\mathcal{L} f_+(x)dx \nonumber
\end{eqnarray}
We only need to treat the two last integrals, and in fact we just need to study one of them since the other can be treated in a similar way. So, for the third integral we have
\begin{eqnarray*}
\int_{\mathbb{R}^n}f_+(x)^{p-2}f_+(x)\mathcal{L} f_-(x)dx&=&\int_{\mathbb{R}^n}f_+(x)^{p-2}f_+(x)\int_{\mathbb{R}^n}[f_-(x)-f_-(y)]\pi(x-y)dydx\\[5mm]
&=&\int_{\mathbb{R}^n}f_+(x)^{p-2}\int_{\mathbb{R}^n}[f_+(x)f_-(x)-f_+(x)f_-(y)]\pi(x-y)dydx
\end{eqnarray*}
However, since $f_+$ and $f_-$ have disjoint supports we obtain the following estimate: 
$$\int_{\mathbb{R}^n}f_+(x)^{p-2}f_+(x)\mathcal{L} f_-(x)dx=-\int_{\mathbb{R}^n}f_+(x)^{p-2}\int_{\mathbb{R}^n}[f_+(x)f_-(y)]\pi(x-y)dydx\leq 0$$
Recalling that $\pi$ is a positive function we obtain that this quantity is negative as all the terms inside the integral are positive. With this observation we see that the last terms of (\ref{FormDecomp}) are positive and we have
\begin{eqnarray*}
\int_{\mathbb{R}^n}f_+(x)^{p-2}f_+(x)\mathcal{L} f_+(x)dx+\int_{\mathbb{R}^n}f_-(x)^{p-2}f_-(x)\mathcal{L} f_-(x)dx \leq \int_{\mathbb{R}^n}|f(x)|^{p-2}f(x)\mathcal{L}f(x)dx<+\infty
\end{eqnarray*}
Then, using the first part of the proof we have $f_\pm\in \dot{B}^{2\alpha/p,p}_{p}(\mathbb{R}^n)$ and since $f=f_+-f_-$ we conclude that $f$ belongs to the Besov space $\dot{B}^{2\alpha/p,p}_{p}(\mathbb{R}^n)$.\hfill$\blacksquare$\\

\textit{\textbf{Proof of Theorem \ref{Theo1}.}} We have obtained with the previous results a family of regular functions $(\theta^{(\varepsilon)})_{\varepsilon >0}\in L^{\infty}([0,T]; L^{p}(\mathbb{R}^n))$ which are solutions of (\ref{SistApprox}) and satisfy the uniform bound $\|\theta^{(\varepsilon)}(\cdot, t)\|_{L^p}\leq \|\theta_0\|_{L^p}$.\\ 
Since $L^{\infty}([0,T]; L^{p}(\mathbb{R}^n))= \left(L^{1}([0,T]; L^{q}(\mathbb{R}^n))\right)'$, with $\frac{1}{p}+\frac{1}{q}=1$, we can extract from those solutions $\theta^{(\varepsilon)}$ a subsequence $(\theta_k)_{k\in \mathbb{N}}$ which is $\ast$-weakly convergent to some function $\theta$ in the space $L^{\infty}([0,T]; L^{p}(\mathbb{R}^n))$, which implies convergence in $\mathcal{D}'(\mathbb{R}^+\times\mathbb{R}^n)$. However, this weak convergence is not sufficient to assure the convergence of $(v_\varepsilon\; \theta_k)$ to $v\;\theta$. For this we use the  remarks that follow.\\

First, using remark \ref{Rem_Approx} we can consider a sequence $(v_k)_{k\in\mathbb{N}}$ with $v_k$ as in formula (\ref{Forbmoaprox}) such that $v_k \longrightarrow v$ weakly in $bmo(\mathbb{R}^n)$. Secondly, combining  Proposition \ref{PropoViscosityMaxPrinc} and Theorem \ref{TheoBesov} we obtain that solutions $\theta_k$ belongs to the space $L^{\infty}([0,T]; L^{p}(\mathbb{R}^n))\cap L^{1}([0,T]; \dot{B}^{2\alpha/p,p}_p(\mathbb{R}^n))$ for all $k\in \mathbb{N}$. \\

To finish, fix a function $\varphi\in \mathcal{C}^{\infty}_{0}([0,T]\times \mathbb{R}^n)$. Then we have the fact that $\varphi \theta_k\in  L^{1}([0,T]; \dot{B}^{2\alpha/p,p}_p(\mathbb{R}^n))$ and $\partial_t \varphi \theta_k\in  L^{1}([0,T]; \dot{B}^{-N,p}_p(\mathbb{R}^n))$. This implies the local inclusion, in space as well as in time, $\varphi \theta_k\in \dot{W}^{2\alpha/p,p}_{t,x}\subset \dot{W}^{2\alpha/p,2}_{t,x}$ so we can apply classical results such as the Rellich's theorem to obtain convergence of $v_k\; \theta_k$ to $v\;\theta$. \\

Thus, we obtain existence and uniqueness of weak solutions for the problem (\ref{Equation0}) with an initial data in $\theta_0\in L^p(\mathbb{R}^n)$, $2\leq p<+\infty$ that satisfy the maximum principle. Moreover, we have that these solutions $\theta(x,t)$ belong to the space $L^{\infty}([0,T]; L^{p}(\mathbb{R}^n))\cap L^{p}([0,T]; \dot{B}^{2\alpha/p,p}_p(\mathbb{R}^n))$.\hfill$\blacksquare$
\begin{Remarque}
These lines explain how to obtain weak solutions from viscosity ones and it will be used freely in the sequel.  
\end{Remarque}
\section{Positivity principle}\label{Sect_PrincipeMax}
We prove in this section Theorem \ref{Theo2}. Recall that by hypothesis we have $0\leq \psi_0\leq M$ an initial datum for the equation (\ref{Equation0}) with $\psi_0\in L^p(\mathbb{R}^n)$ and $\frac{n}{2\min(\beta, \delta)}\leq p\leq+\infty$.\\

To begin with, we fix two constants, $\rho, R$ such that $R>2\rho>0$. Then we set $A_{0,R}(x)$ a function equals to $M/2$ over $|x|\leq 2R$ and equals to $\psi_0(x)$ over $|x|>2R$ and we write $B_{0,R}(x)=\psi_0(x)-A_{0,R}(x)$, so by construction we have $$\psi_0(x)=A_{0,R}(x)+B_{0,R}(x)$$ with $\|A_{0,R}\|_{L^\infty} \leq M$ and $\|B_{0,R}\|_{L^\infty} \leq M/2$. Remark that $A_{0,R}, B_{0,R}\in L^p(\mathbb{R}^n)$.\\

Now fix $v\in L^{\infty}([0,T];bmo(\mathbb{R}^n))$ such that $div(v)=0$ and consider the equations
\begin{equation}\label{ForDouble}
\begin{cases}
\partial_t A_R(x,t)+ \nabla\cdot (v\,A_R)(x,t)+\mathcal{L} A_R(x,t)=0,\\[4mm]
A_R(x,0)=A_{0,R}(x)
\end{cases}
\mbox{and}\qquad
\begin{cases}
\partial_t B_R(x,t)+ \nabla\cdot (v\,B_R)(x,t)+\mathcal{L} B_R(x,t)=0\\[4mm]
B_R(x,0)= B_{0,R}(x).
\end{cases}
\end{equation}
Using the maximum principle and by construction we have the following estimates for $t\in [0,T]$:
\begin{eqnarray}
\|A_R(\cdot, t)\|_{L^p}&\leq& \|A_{0,R}\|_{L^p}\leq \|\psi_0\|_{L^p}+CM R^{n/p} \quad (1<p<+\infty)\label{Formula1}\\[5mm] 
\|A_R(\cdot, t)\|_{L^\infty}&\leq & \|A_{0,R}\|_{L^\infty}\leq M.\nonumber\\[5mm]
\|B_R(\cdot, t)\|_{L^\infty}&\leq & \|B_{0,R}\|_{L^\infty}\leq M/2. \nonumber
\end{eqnarray}
where $A_R(x,t)$ and $B_R(x,t)$ are the weak solutions of the systems (\ref{ForDouble}).
\begin{Lemme}\label{lemmeRecollement} 
The function $\psi(x,t)=A_R(x,t)+B_R(x,t)$ is the unique solution for the problem
\begin{equation}\label{Equation1}
\left\lbrace
\begin{array}{l}
\partial_t \psi(x,t)+ \nabla\cdot (v\,\psi)(x,t)+\mathcal{L} \psi(x,t)=0\\[5mm]
\psi(x,0)=A_{0,R}(x)+B_{0,R}(x).
\end{array}
\right.
\end{equation}
\end{Lemme}
\textit{\textbf{Proof}.} Using hypothesis for $A_R(x,t)$ and $B_R(x,t)$ and the linearity of equation (\ref{Equation1}) we have that the function $\psi_R(x,t)=A_R(x,t)+B_R(x,t)$ is a solution for this equation. Uniqueness is assured by the maximum principle and by the continuous dependence from initial data given in corollary \ref{CorDepContinue}, thus we can write $\psi(x,t)=\psi_R(x,t)$.\hfill $\blacksquare$\\

To continue, we will need an auxiliary function $\phi \in \mathcal{C}^{\infty}_{0}(\mathbb{R}^n)$ such that $\phi(x)=0$ for $|x|\geq 1$ and $\phi(x)=1$ if $|x|\leq 1/2$ and we set $\varphi(x)=\phi(x/R)$. Now, we will estimate the $L^p$-norm of $\varphi(x)(A_R(x,t)-M/2)$ with $p>n/2\min(\beta,\delta)$, where $\beta$ and $\delta$ are the parameters of the hypothesis for the function $\pi$ in the cases \textbf{(a)-(d)}. We write:
\begin{eqnarray}
\partial_t\|\varphi(\cdot)(A_R(\cdot,t)-M/2)\|_{L^p}^p &=& p\int_{\mathbb{R}^n}\big|\varphi(x)(A_R(x,t)-M/2)\big|^{p-2}\big(\varphi(x)(A_R(x,t)-M/2) \big)\nonumber\\[4mm]
& & \times\; \partial_t\big(\varphi(x)(A_R(x,t)-M/2) \big) dx \label{equat1}
\end{eqnarray}
We observe that we have the following identity for the last term above
\begin{eqnarray*}
\partial_t(\varphi(x)(A_R(x,t)-M/2))&=&- \nabla \cdot (\varphi(x) \, v(A_R(x,t)-M/2))-\mathcal{L} (\varphi(x)(A_R(x,t)-M/2))\\[3mm]
&+&(A_R(x,t)-M/2)v\cdot \nabla \varphi(x)+[\mathcal{L}, \varphi]A_R(x,t)-M/2 \mathcal{L} \varphi(x)
\end{eqnarray*}
where we noted $[\mathcal{L}, \varphi]$ the commutator between $\mathcal{L}$ and $\varphi$. Thus, using this identity in (\ref{equat1}) and the fact that $div(v)=0$ we have
\begin{eqnarray}
\partial_t\|\varphi(\cdot)(A_R(\cdot,t)-M/2)\|_{L^p}^p &=&-p\int_{\mathbb{R}^n}\big|\varphi(x)(A_R(x,t)-M/2)\big|^{p-2}\big(\varphi(x)(A_R(x,t)-M/2) \big)\nonumber\\[4mm]
& & \times\; \mathcal{L}(\varphi(x)(A_R(x,t)-M/2))dx\label{equa2}\\[4mm]
&+&p\int_{\mathbb{R}^n}\big|\varphi(x)(A_R(x,t)-M/2)\big|^{p-2}\big(\varphi(x)(A_R(x,t)-M/2) \big)\nonumber\\[4mm]
& & \times\; \left([\mathcal{L}, \varphi]A_R(x,t)-M/2\mathcal{L} \varphi(x)\right)dx\nonumber
\end{eqnarray}
Remark that the integral (\ref{equa2}) is positive so one has
\begin{eqnarray*}
\partial_t\|\varphi(\cdot)(A_R(\cdot,t)-M/2)\|_{L^p}^p &\leq &p\int_{\mathbb{R}^n}\big|\varphi(x)(A_R(x,t)-M/2)\big|^{p-2}\big(\varphi(x)(A_R(x,t)-M/2) \big)\\[4mm]
& & \times\; \left([\mathcal{L}, \varphi]A_R(x,t)-M/2 \mathcal{L}\varphi(x)\right)dx
\end{eqnarray*}
Using Hölder's inequality and integrating in time the previous expression we have
\begin{eqnarray}\label{EquaPrevious}
\|\varphi(\cdot)(A_R(\cdot,t)-M/2)\|^p_{L^p}& \leq &\|\varphi(\cdot)(A_R(\cdot,0)-M/2)\|^p_{L^p}+\int_{0}^t\bigg(\left\|[\mathcal{L}, \varphi]A_R(\cdot,s)\right\|_{L^p}+ \|M/2 \mathcal{L}\varphi\|_{L^p}\bigg)ds
\end{eqnarray}
The first term of the right side is null since over the support of $\varphi$ we have identity $A_R(x,0)=M/2$. For the second term $\left\|[\mathcal{L}, \varphi]A_R(\cdot,s)\right\|_{L^p}$ we will need the following lemma
\begin{Lemme}
For $1\leq p\leq +\infty$ we have for the cases \textbf{(a)-(d)} the following inequality:
\begin{equation*}
\left\|[\mathcal{L}, \varphi]A_R(\cdot,s)\right\|_{L^p}\leq C(R^{-2\beta}+R^{-2\delta})\|A_{0,R}\|_{L^p}.
\end{equation*}
\end{Lemme}
\textbf{\textit{Proof.}} We have $[\mathcal{L}, \varphi]A_R(x,s)=\displaystyle{\int_{\mathbb{R}^n}}\big(\varphi(x)-\varphi(x-y)\big)A_R(x-y,s)\pi(y)dy$ and we divide our study following the different cases \textbf{(a)-(d)}.\\

For the case \textbf{(a)}, where $0<\alpha\leq \beta<1/2$ and $0<\delta<1/2$, or in the case \textbf{(b)} where $0<\alpha=\beta=\delta<1/2$, we proceed as follows. We begin with the case $p=+\infty$ and we write:
\begin{equation}\label{Equation103}
|[\mathcal{L}, \varphi]A_R(x,s)|\leq \int_{\mathbb{R}^n}\frac{|\varphi(x)-\varphi(y)|}{|x-y|^{n+2\beta}}|A_R(y,s)|dy+\int_{\mathbb{R}^n}\frac{|\varphi(x)-\varphi(y)|}{|x-y|^{n+2\delta}}|A_R(y,s)|dy
\end{equation}
Again, it is enough to study one of these two integrals since the other can be treated in a totally similar way. We write:
\begin{eqnarray*}
\int_{\mathbb{R}^n}\frac{|\varphi(x)-\varphi(y)|}{|x-y|^{n+2\beta}}|A_R(y,s)|dy &=& \int_{\{|x-y|>R\}}\frac{|\varphi(x)-\varphi(y)|}{|x-y|^{n+2\beta}}|A_R(y,s)|dy+\int_{\{|x-y|\leq R\}}\frac{|\varphi(x)-\varphi(y)|}{|x-y|^{n+2\beta}}|A_R(y,s)|dy\\
&\leq & 2\|\varphi\|_{L^\infty} \int_{\{|x-y|>R\}}\frac{|A_R(y,s)|}{|x-y|^{n+2\beta}}dy+\int_{\{|x-y|\leq R\}}\frac{\| \nabla \varphi\|_{L^\infty}|x-y|}{|x-y|^{n+2\beta}}|A_R(y,s)|dy\\
&\leq & 2\|\varphi\|_{L^\infty}\|A_R(\cdot,s)\|_{L^\infty} \int_{\{|x-y|>R\}}\frac{1}{|x-y|^{n+2\beta}}dy+ C R^{-1}\int_{\{|x-y|\leq R\}}\frac{|A_R(y,s)|}{|x-y|^{n+2\beta-1}}dy\\
&\leq & 2C\|\varphi\|_{L^\infty}\|A_R(\cdot,s)\|_{L^\infty} R^{-2\beta}+C \|A_R(\cdot,s)\|_{L^\infty} R^{-2\beta}\leq CR^{-2\beta}\|A_{0,R}\|_{L^\infty}. 
\end{eqnarray*}
Then, with the $\delta$-part in inequality (\ref{Equation103}) we have
$$\|[\mathcal{L}, \varphi]A_R(\cdot,s)\|_{L^\infty}\leq C(R^{-2\beta}+R^{-2\delta})\|A_{0,R}\|_{L^\infty}.$$
The case $p=1$ is very similar. Using inequality (\ref{Equation103}) we have
$$\int_{\mathbb{R}^n}|[\mathcal{L}, \varphi]A_R(x,s)|dx\leq \int_{\mathbb{R}^n}\int_{\mathbb{R}^n}\frac{|\varphi(x)-\varphi(y)|}{|x-y|^{n+2\beta}}|A_R(y,s)|dydx+\int_{\mathbb{R}^n}\int_{\mathbb{R}^n}\frac{|\varphi(x)-\varphi(y)|}{|x-y|^{n+2\delta}}|A_R(y,s)|dydx
$$
We only estimate one of the previous integrals.
\begin{eqnarray*}
\int_{\mathbb{R}^n}\int_{\mathbb{R}^n}\frac{|\varphi(x)-\varphi(y)|}{|x-y|^{n+2\beta}}|A_R(y,s)|dydx &\leq & C\|\varphi\|_{L^\infty} \int_{\mathbb{R}^n}\int_{\{|x-y|>R\}}\frac{|A_R(y,s)|}{|x-y|^{n+2\beta}}dydx\\
& &+ R^{-1}\int_{\mathbb{R}^n}\int_{\{|x-y|\leq R\}}\frac{|A_R(y,s)|}{|x-y|^{n+2\beta-1}}dydx\\
&\leq & C\|\varphi\|_{L^\infty} \|A_R(\cdot,s)\|_{L^1}R^{-2\beta}+C\|A_R(\cdot,s)\|_{L^1}R^{-2\beta}\leq CR^{-2\beta}\|A_{0,R}\|_{L^1}.
\end{eqnarray*}
With the other integral, we obtain
$$\|[\mathcal{L}, \varphi]A_R(\cdot,s)\|_{L^1}\leq C(R^{-2\beta}+R^{-2\delta})\|A_{0,R}\|_{L^1}.$$
Finally, the case $1<p<+\infty$ is obtained by interpolation. See \cite{Grafakos} or \cite{Stein2} for more details about interpolation.\\

For the remaining cases \textbf{(c)} and \textbf{(d)} (\textit{i.e.} if $\alpha=\beta=1/2$ and $0<\delta<1/2$ or $\alpha=\beta=\delta=1/2$), the result will be a consequence of the Calder\'on's commutator inequality (see \cite{Grafakos}) and the maximum principle. \hfill$\blacksquare$\\

Now, getting back to the last term of (\ref{EquaPrevious}) we have by the definition of $\varphi$ and the properties of the operator $\mathcal{L}$ the estimate:
$$\|M/2 \mathcal{L}\varphi\|_{L^p}\leq C MR^{n/p}(R^{-2\beta}+R^{-2\delta}).$$
We thus have
$$\|\varphi(\cdot)(A_R(\cdot,t)-M/2)\|^p_{L^p}\leq C(R^{-2\beta}+R^{-2\delta})\int_{0}^t\bigg(\|A_{0,R}\|_{L^p}+MR^{n/p}\bigg)ds.$$
Observe that we have at our disposal estimate (\ref{Formula1}), so we can write
$$\|\varphi(\cdot)(A_R(\cdot,t)-M/2)\|^p_{L^p}\leq Ct (R^{-2\beta}+R^{-2\delta})\left(\|\psi_0\|_{L^p}+MR^{n/p}\right)$$
Using again the definition of $\varphi$ one has $\displaystyle{\int_{B(0,\rho)}}|A_R(\cdot,t)-M/2|^{p}dx\leq Ct(R^{-2\beta}+R^{-2\delta})\left(\|\psi_0\|_{L^p}+MR^{n/p}\right)$. Thus, if $R\longrightarrow +\infty$ and since $p>\frac{n}{2\min(\beta,\delta)}$, we have $A(x,t)=M/2$ over $B(0,\rho)$.\\ 

Hence, by construction we have  $\psi(x,t)=A_R(x,t)+B_R(x,t)$ where $\psi$ is a solution of (\ref{Equation1}) with initial data $\psi_0=A_{0,R}+B_{0,R}$, but, since over $B(0,\rho)$ we have $A(x,t)=M/2$ and $\|B(\cdot,t)\|_{L^\infty}\leq M/2$, one finally has the desired estimate $0\leq \psi(x,t)\leq M$.\hfill$\blacksquare$
\section{Existence of solutions with a $L^\infty$ initial data}\label{SecLinfty}

The proof given before for the positivity principle allows us to obtain the existence of solutions for the fractional diffusion transport equation (\ref{Equation0}) when the initial data $\theta_0$ belongs to the space $L^\infty(\mathbb{R}^n)$. The utility of this fact will appear clearly in the next section as it will be used in Theorem \ref{Theo3}.\\

Let us fix $\theta_0^R=\theta_0 \mathds{1}_{B(0,R)}$ with $R>0$ so we have $\theta_0^R\in L^p(\mathbb{R}^n)$ for all $1\leq p\leq +\infty$. Following section \ref{SecExiUnic}, there is a unique solution $\theta^R$ for the problem
\begin{equation*}
\left\lbrace
\begin{array}{l}
\partial_t \theta^R+ \nabla\cdot (v\theta^R)+\mathcal{L}\theta^R=0\\[5mm]
\theta^R(x,0)=\theta_0^R(x)\\[5mm]
div(v)=0 \quad \mbox{ and } v\in L^{\infty}([0,T];  bmo(\mathbb{R}^n)).
\end{array}
\right.
\end{equation*}
such that $\theta^R\in L^\infty([0,T]; L^p(\mathbb{R}^n))$. By the maximum principle we have $\|\theta^R(\cdot, t)\|_{L^p}\leq\|\theta^R_0\|_{L^p}\leq v_n\|\theta_0\|_{L^\infty}R^{n/p}$. Taking the limit $p\longrightarrow +\infty$ and making $R\longrightarrow +\infty$ we finally get
$$\|\theta(\cdot, t)\|_{L^\infty}\leq C \|\theta_0\|_{L^\infty}.$$
This shows that for an initial data $\theta_0\in  L^\infty(\mathbb{R}^n)$ there exists an associated solution $\theta\in L^\infty([0,T];L^\infty(\mathbb{R}^n))$.
\section{H\"older Regularity}\label{SeccHolderRegularity}
In this section we are going to prove Theorem \ref{Theo3}. It is very important to note that we will work only with the cases \textbf{(c)} and \textbf{(d)}: from now on the operator $\mathcal{L}$ is assumed to be of the form (\ref{DefOperator1}) with an associated Lévy measure $\pi$ satisfiying the hypothesis (\ref{DefKernel2}) and (\ref{DefKernel3}) with $\alpha=\beta=1/2$ and $0<\delta<1/2$ or $\alpha=\beta=\delta=1/2$.\\

We will now study H\"older regularity by duality using Hardy spaces. These spaces have several equivalent characterizations (see \cite{Coifmann}, \cite{Gold} and \cite{Stein2} for a detailed treatment). In this paper we are interested mainly in the molecular approach that defines \textit{local} Hardy spaces.
\begin{Definition}[Local Hardy spaces $h^\sigma$] Let $0<\sigma<1$. The local Hardy space $h^{\sigma}(\mathbb{R}^n)$ is the set of distributions $f$ that admits the following molecular decomposition:
\begin{equation}\label{MolDecomp}
f=\sum_{j\in \mathbb{N}}\lambda_j \psi_j
\end{equation}
where $(\lambda_j)_{j\in \mathbb{N}}$ is a sequence of complex numbers such that $\sum_{j\in \mathbb{N}}|\lambda_j|^\sigma<+\infty$ and $(\psi_j)_{j\in \mathbb{N}}$ is a family of $r$-molecules in the sense of definition \ref{DefMolecules} below. The $h^\sigma$-norm\footnote{it is not actually a \textit{norm} since $0<\sigma<1$. More details can be found in \cite{Gold} and \cite{Stein2}.} is then fixed by the formula 
$$\|f\|_{h^\sigma}=\inf\left\{\left(\sum_{j\in \mathbb{N}}|\lambda_j|^\sigma\right)^{1/\sigma}:\; f=\sum_{j\in \mathbb{N}}\lambda_j \psi_j \right\}$$ where the infimum runs over all possible decompositions (\ref{MolDecomp}).
\end{Definition}
Local Hardy spaces have many remarquable properties and we will only stress here, before passing to duality results concerning $h^\sigma$ spaces, the fact that Schwartz class $\mathcal{S}(\mathbb{R}^n)$ is dense in $h^{\sigma}(\mathbb{R}^n)$. \\

Now, let us take a closer look at the dual space of the local Hardy spaces. In \cite{Gold} D. Goldberg proved the following important theorem:
\begin{Theoreme}[Hardy-Hölder duality]\label{TheoHHD} Let $\frac{n}{n+1}<\sigma<1$ and fix $\gamma=n(\frac{1}{\sigma}-1)$. Then the dual of local Hardy space $h^{\sigma}(\mathbb{R}^n)$ is the Hölder space $\mathcal{C}^\gamma(\mathbb{R}^n)$ fixed by the norm
$$\|f\|_{\mathcal{C}^\gamma}=\|f\|_{L^\infty}+\underset{x\neq y}{\sup}\frac{|f(x)-f(y)|}{|x-y|^\gamma}.$$
\end{Theoreme}
This result allows us to study the Hölder regularity of functions in terms of Hardy spaces and it will be applied to the solutions of the equation (\ref{Equation0}).
\begin{Remarque}\label{Remark3}
Since $\frac{n}{n+1}<\sigma<1$, we have $\sum_{j\in \mathbb{N}}|\lambda_j|\leq\left(\sum_{j\in \mathbb{N}}|\lambda_j|^\sigma\right)^{1/\sigma}$ thus for testing Hölder continuity of a function $f$ it is enough to study the quantities $|\langle f,\psi_j\rangle|$ where $\psi_j$ is an $r$-molecule.
\end{Remarque}
Since we are going to work with local Hardy spaces, we will introduce a size treshold in order to distinguish \textit{small} molecules from \textit{big} ones in the following way:
\begin{Definition}[$r$-molecules]\label{DefMolecules} Set $\frac{n}{n+1}<\sigma<1$, define $\gamma=n(\frac{1}{\sigma}-1)$ and fix a real number $\omega$ such that $0<\gamma<\omega<1$. An integrable function $\psi$ is an $r$-molecule if we have
\begin{enumerate}
\item[$\bullet$]\underline{Small molecules $(0<r<1)$:}
\begin{eqnarray}
& & \int_{\mathbb{R}^n} |\psi(x)||x-x_0|^{\omega}dx \leq  r^{\omega-\gamma}\mbox{, for } x_0\in \mathbb{R}^n\label{Hipo1} \;\qquad\qquad\qquad\mbox{(concentration condition)} \\[5mm]
& &\|\psi\|_{L^\infty}  \leq  \frac{1}{r^{n+\gamma}}\label{Hipo2} \qquad\qquad\qquad\qquad\qquad\qquad\qquad\qquad\qquad\, \mbox{(height condition)}  \\[5mm]
& &\int_{\mathbb{R}^n} \psi(x)dx=0\label{Hipo3}\qquad\qquad\qquad\qquad\qquad\qquad\qquad\qquad\qquad \mbox{(moment condition)} 
\end{eqnarray}
\item[$\bullet$] \underline{Big molecules $(1\leq r<+\infty)$:}\\

In this case we only require conditions (\ref{Hipo1}) and (\ref{Hipo2}) for the $r$-molecule $\psi$ while the moment condition (\ref{Hipo3}) is dropped.
\end{enumerate}
\end{Definition}

\begin{Remarque}\label{Remark2}
\begin{itemize}
\item[]
\item[1)] Note that the point $x_0\in \mathbb{R}^n$ can be considered as the ``center'' of the molecule.
\item[2)] Conditions (\ref{Hipo1}) and (\ref{Hipo2}) imply the estimate $\|\psi\|_{L^1}\leq C\, r^{-\gamma}$ thus every $r$-molecule belongs to $L^p(\mathbb{R}^n)$ with $1<p<+\infty$.
\item[3)] In this definition, we find convenient to show explicitely the dependence on the H\"older parameter $\gamma$ instead of $\sigma$. 
\end{itemize}
\end{Remarque}
The main interest of using molecules relies in the possibility of \textit{transfering} the regularity problem to the evolution of such molecules:

\begin{Proposition}[Transfer property]\label{Transfert} Let $\psi(x,s)$ be a solution of the backward problem
\begin{equation}
\left\lbrace
\begin{array}{rl}
\partial_s \psi(x,s)=& -\nabla\cdot [v(x,t-s)\psi(x,s)]-\mathcal{L}\psi(x,s)\label{Evolution1}\\[5mm]
\psi(x,0)=& \psi_0(x)\in L^1\cap L^\infty(\mathbb{R}^n)\\[5mm]
div(v)=&0 \quad \mbox{and }\; v\in L^{\infty}([0,T];bmo(\mathbb{R}^n))
\end{array}
\right.
\end{equation}
If $\theta(x,t)$ is a solution of (\ref{Equation0}) with $\theta_0\in L^\infty(\mathbb{R}^n)$ then we have the identity
\begin{equation*}
\int_{\mathbb{R}^n}\theta(x,t)\psi(x,0)dx=\int_{\mathbb{R}^n}\theta(x,0)\psi(x,t)dx.
\end{equation*}
\end{Proposition}
\textit{\textbf{Proof.}}
We first consider the expression
$$\partial_s\int_{\mathbb{R}^n}\theta(x,t-s)\psi(x,s)dx=\int_{\mathbb{R}^n}-\partial_s\theta(x,t-s)\psi(x,s)+\partial_s\psi(x,s)\theta(x,t-s)dx.$$
Using equations (\ref{Equation0}) and (\ref{Evolution1}) we obtain
\begin{eqnarray*}
\partial_s\int_{\mathbb{R}^n}\theta(x,t-s)\psi(x,s)dx&=&\int_{\mathbb{R}^n}- \nabla\cdot\left[(v(x,t-s)\theta(x,t-s)\right]\psi(x,s)+\mathcal{L}\theta(x,t-s)\psi(x,s) \\[5mm]
&-& \nabla\cdot\left[(v(x,t-s) \psi(x,s))\right]\theta(x,t-s)-\mathcal{L}\psi(x,s) \theta(x,t-s)dx.
\end{eqnarray*}
Now, using the fact that $v$ is divergence free and the symmetry of the operator $\mathcal{L}$ we have that the expression above is equal to zero, so the quantity
$$\int_{\mathbb{R}^n}\theta(x,t-s)\psi(x,s)dx$$
remains constant in time. We only have to set $s=0$ and $s=t$ to conclude. \hfill$\blacksquare$\\

This proposition says, that in order to control $\langle \theta(\cdot, t),\psi_0\rangle$, it is enough (and much simpler) to study the bracket $\langle \theta_0,\psi(\cdot, t)\rangle$. \\

\textbf{Proof of Theorem \ref{Theo3}.} 
Once we have the transfer property proven above, the proof of Theorem \ref{Theo3} is quite direct and it reduces to a $L^1$ estimate for molecules. Indeed, assume that for \textit{all} molecular initial data $\psi_0$ we have a $L^1$ control for $\psi(\cdot,t)$ a solution of (\ref{Evolution1}), then Theorem \ref{Theo3} follows easily: applying Proposition \ref{Transfert} with the fact that $\theta_0\in L^{\infty}(\mathbb{R}^n)$ we have 
\begin{equation}\label{DualQuantity1}
|\langle \theta(\cdot, t), \psi_0\rangle|=\left|\int_{\mathbb{R}^n}\theta(x,t)\psi_0(x)dx\right|=\left|\int_{\mathbb{R}^n}\theta(x,0)\psi(x,t)dx\right|\leq \|\theta_0\|_{L^\infty}\|\psi(\cdot,t)\|_{L^1}<+\infty.
\end{equation}
From this, we obtain that $\theta(\cdot, t)$ belongs to the Hölder space $\mathcal{C}^\gamma(\mathbb{R}^n)$.\\

Now we need to study the control of the $L^1$ norm of $\psi(\cdot, t)$ and we divide our proof in two steps following the molecule's size. For the initial big molecules, \textit{i.e.} if $r\geq 1$, the needed control is straightforward: apply the maximum principle and the remark \ref{Remark2}-2) above to obtain
$$\|\theta_0\|_{L^\infty}\|\psi(\cdot,t)\|_{L^1}\leq \|\theta_0\|_{L^\infty}\|\psi_0\|_{L^1}\leq C \frac{1}{r^\gamma} \|\theta_0\|_{L^\infty},$$
but, since $r\geq 1$, we have that $|\langle \theta(\cdot, t), \psi_0\rangle|<+\infty$ for all \textit{big} molecules.\\

In order to finish the proof of this theorem, it  only remains to treat the $L^1$ control for \textit{small} molecules. This is the most complex part of the proof and it is treated in the following theorem:
\begin{Theoreme}\label{TheoL1control}
For all small $r$-molecules (\textit{i.e. }$0<r<1$), there exists a time $T_0>0$ such that we have the following control of the $L^1$-norm.
$$\|\psi(\cdot, t)\|_{L^1}\leq C T_0^{-\gamma}\qquad (T_0<t<T),$$
with $0<\gamma<2\delta\leq 1$. 
\end{Theoreme}
This theorem will be proven in sections \ref{SecEvolMol1}, \ref{SecEvolMol2} and \ref{SecIterationMol}.\\ 

Accepting for a while this result, we have then a good control over the quantity $\|\psi(\cdot, t)\|_{L^1}$  for all $0<r<1$ and getting back to (\ref{DualQuantity1}) we obtain that $|\langle \theta(\cdot, t), \psi_0\rangle|$ is always bounded for $T_0<t<T$ and for any molecule $\psi_0$: we have proven by a duality argument the Theorem \ref{Theo3}. \hfill $\blacksquare$\\

Let us now briefly explain the main steps of Theorem \ref{TheoL1control}. We need to construct a suitable control in time for the $L^1$-norm of the solutions $\psi(\cdot,t)$ of the backward problem (\ref{Evolution1}) where the inital data $\psi_0$ is a \textit{small} $r$-molecule. This will be achieved by iteration in two different steps. The first step explains the molecules' deformation after a very small time $s_0>0$, which is related to the size $r$ by the bounds $0<s_0\leq \epsilon r$ with $\epsilon$ a small constant. In order to obtain a control of the $L^1$ norm for larger times we need to perform a second step which takes as a starting point the results of the first step and gives us the deformation for another small time $s_1$, which is also related to the original size $r$. Once this is achieved it is enough to iterate the second step as many times as necessary to get rid of the dependence of the times $s_0, s_1,...$ from the molecule's size. This way we obtain the $L^1$ control needed for all time $T_0<t<T$.\\

\subsection{Small time molecule's evolution: First step}\label{SecEvolMol1}
The following theorem shows how the molecular properties are deformed with the evolution for a small time $s_0$.
\begin{Theoreme}\label{SmallGeneralisacion} Set $\sigma$, $\gamma$ and $\omega$ three real numbers such that $\frac{n}{n+1}<\sigma<1$, $\gamma=n(\frac{1}{\sigma}-1)$ and $0<\gamma<\omega<2\delta<1$ in the case \textbf{(c)} or $0<\gamma<\omega<1$ in the case \textbf{(d)}. Let $\psi(x,s_0)$ be a solution of the problem
\begin{equation}\label{SmallEvolution}
\left\lbrace
\begin{array}{rl}
\partial_{s_0} \psi(x,s_0)=& -\nabla\cdot(v\, \psi)(x,s_0)-\mathcal{L}\psi(x,s_0)\\[5mm]
\psi(x,0)=& \psi_0(x)\\[5mm]
div(v)=&0 \quad \mbox{and }\; v\in L^{\infty}([0,T];bmo(\mathbb{R}^n))\quad \mbox{with } \underset{s_0\in [0,T]}{\sup}\; \|v(\cdot,s_0)\|_{bmo}\leq \mu
\end{array}
\right.
\end{equation}
If $\psi_0$ is a small $r$-molecule in the sense of definition \ref{DefMolecules} for the local Hardy space $h^\sigma(\mathbb{R}^n)$, then there exists a positive constant $K=K(\mu)$ big enough and a positive constant $\epsilon$ such that for all $0< s_0 \leq\epsilon r$ small we have the following estimates
\begin{eqnarray}
\int_{\mathbb{R}^n}|\psi(x,s_0)||x-x(s_0)|^\omega dx &\leq &(r+Ks_0)^{\omega-\gamma}  \label{SmallConcentration}\\
\|\psi(\cdot, s_0)\|_{L^\infty}&\leq & \frac{1}{\big(r+Ks_0\big)^{n+\gamma}}\label{SmallLinftyevolution}\\
\|\psi(\cdot, s_0)\|_{L^1} &\leq & \frac{v_n}{\big(r+Ks_0\big)^{\gamma}}\label{SmallL1evolution}
\end{eqnarray}
where $v_n$ denotes the volume of the $n$-dimensional unit ball.\\ 

The new molecule's center $x(s_0)$ used in formula (\ref{SmallConcentration}) is fixed by 
\begin{equation}\label{Defpointx_0}
\left\lbrace
\begin{array}{rl}
x'(s_0)=& \overline{v}_{B_r}=\frac{1}{|B_r|}\displaystyle{\int_{B_r}}v(y,s_0)dy \quad \mbox{where } B_r=B(x(s_0),r).\\[5mm]
x(0)=& x_0.
\end{array}
\right.
\end{equation}  
\end{Theoreme}

\begin{Remarque}
\begin{itemize}
\item[]
\item[1)] The definition of the point $x(s_0)$ given by (\ref{Defpointx_0}) reflects the molecule's center transport using velocity $v$.
\item[2)] Remark that it is enough to treat the case $0<(r+Ks_0)<1$ since $s_0$ is small: otherwise the $L^1$ control will be trivial for time $s_0$ and beyond: we only need to apply the maximum principle.
\end{itemize}
\end{Remarque}

The proof of this theorem follows the next scheme: the small concentration condition (\ref{SmallConcentration}), which is proven in the Proposition \ref{Propo1}, implies the height condition (\ref{SmallLinftyevolution}) (proved in Proposition \ref{Propo2}). Once we have these two  conditions, the $L^1$ estimate (\ref{SmallL1evolution}) will follow easily and this is proven in Proposition \ref{Propo3}.
\begin{Proposition}[Small time Concentration condition]\label{Propo1}
Under the  hypothesis of Theorem \ref{SmallGeneralisacion}, if $\psi_0$ is a small $r$-molecule, then the solution $\psi(x,s)$ of (\ref{SmallEvolution}) satisfies
\begin{equation*}
\int_{\mathbb{R}^n} |\psi(x,s_0)||x-x(s_0)|^{\omega}dx \leq (r+Ks_0)^{\omega-\gamma}
\end{equation*}
for $x(s_0)\in \mathbb{R}^n$ fixed by formula (\ref{Defpointx_0}) and with $0<s_0\leq \epsilon r$.
\end{Proposition}
\textit{\textbf{Proof}.}
Let us write $\Omega_0(x)=|x-x(s_0)|^{\omega}$ and $\psi(x)=\psi_+(x)-\psi_-(x)$ where the functions $\psi_{\pm}(x)\geq 0$ have disjoint support. We will note $\psi_\pm(x,s_0)$ solutions of (\ref{SmallEvolution}) with $\psi_\pm(x,0)=\psi_\pm(x)$. At this point, we use the positivity principle, thus by linearity we have
$$|\psi(x,s_0)|=|\psi_+(x,s_0)-\psi_-(x,s_0)|\leq \psi_+(x,s_0)+\psi_-(x,s_0)$$ and we can write
$$\int_{\mathbb{R}^n}|\psi(x,s_0)|\Omega_0(x)dx\leq\int_{\mathbb{R}^n}\psi_+(x,s_0)\Omega_0(x)dx+\int_{\mathbb{R}^n}\psi_-(x,s_0)\Omega_0(x)dx$$
so we only have to treat one of the integrals on the right side above. We have:
\begin{eqnarray*}
I&=&\left|\partial_{s_0} \int_{\mathbb{R}^n}\Omega_0(x)\psi_+(x,s_0)dx\right|\\
&=&\left|\int_{\mathbb{R}^n}\partial_{s_0} \Omega_0(x)\psi_+(x,s_0)+\Omega_0(x)\left[-\nabla\cdot(v\, \psi_+(x,s_0))-\mathcal{L}\psi_+(x,s_0)\right]dx\right|\\
&=&\left|\int_{\mathbb{R}^n}-\nabla\Omega_0(x)\cdot x'(s_0)\psi_+(x,s_0)+\Omega_0(x)\left[-\nabla\cdot(v\, \psi_+(x,s_0))-\mathcal{L}\psi_+(x,s_0)\right]dx\right|\\
&
\end{eqnarray*}
Using the fact that $v$ is divergence free, we obtain
\begin{equation*}
I=\left|\int_{\mathbb{R}^n}\nabla\Omega_0(x)\cdot(v-x'(s_0))\psi_+(x,s_0)-\Omega_0(x)\mathcal{L}\psi_+(x,s_0)dx\right|.
\end{equation*}
Since the operator $\mathcal{L}$ is symmetric and using the definition of $x'(s_0)$ given in (\ref{Defpointx_0}) we have
\begin{equation}\label{smallEstrella}
I\leq c \underbrace{\int_{\mathbb{R}^n}|x-x(s_0)|^{\omega-1}|v-\overline{v}_{B_r}| |\psi_+(x,s_0)|dx}_{I_1} + c\underbrace{\int_{\mathbb{R}^n}\big|\mathcal{L}\Omega_0(x)\big|\, |\psi_+(x,s_0)|dx}_{I_2}.
\end{equation}
We will study separately each of the integrals $I_1$ and $I_2$ in the Lemmas \ref{Lemme1} and \ref{Lemme2} below. But before, we will need the following result
\begin{Lemme}\label{PropoBMO1} Let $f\in bmo(\mathbb{R}^n)$, then
\begin{enumerate}
\item[1)] for all $1<p<+\infty$, $f$ is locally in $L^p$ and $\frac{1}{|B|}\displaystyle{\int_{B}|f(x)-f_B|^pdx}\leq C \|f\|_{bmo}^p$
\item[2)] for all $k\in \mathbb{N}$, we have $|f_{2^k B}-f_B|\leq Ck \|f\|_{bmo}$ where $2^kB=B(x,2^k R)$ is a ball centered at a point $x$ of radius $2^kR$. 
\end{enumerate}
\end{Lemme}
For a proof of these results see \cite{Stein2}.
\begin{Lemme}\label{Lemme1} For integral $I_1$ above we have the estimate $I_1\leq C \mu \; r^{\omega-1-\gamma}$.
\end{Lemme}
\textit{\textbf{Proof}.} We begin by considering the space $\mathbb{R}^n$ as the union of a ball with dyadic coronas centered around $x(s_0)$, more precisely we set $\mathbb{R}^n=B_r\cup \bigcup_{k\geq 1}E_k$ where
\begin{equation}\label{SmallDecoupage}
B_r= \{x\in \mathbb{R}^n: |x-x(s_0)|\leq r\} \quad \mbox{and}\quad E_k= \{x\in \mathbb{R}^n: r2^{k-1}<|x-x(s_0)|\leq r 2^{k}\} \quad \mbox{for } k\geq 1,
\end{equation}
\begin{enumerate}
\item[(i)] \underline{Estimations over the ball $B_r$}. Applying Hölder's inequality to the integral $I_{1,B_r}$ we obtain
\begin{eqnarray}
I_{1,B_r}=\int_{B_r}|x-x(s_0)|^{\omega-1}|v-\overline{v}_{B_r}| |\psi_+(x,s_0)|dx & \leq &\underbrace{\||x-x(s_0)|^{\omega-1}\|_{L^p(B_r)}}_{(1)} \label{Equa1} \\
& \times &\underbrace{\|v-\overline{v}_{B_r}\|_{L^z(B_r)}}_{(2)}\underbrace{\|\psi_+(\cdot, s_0)\|_{L^q(B_r)}}_{(3)}\nonumber
\end{eqnarray}
where $\frac{1}{p}+\frac{1}{z}+\frac{1}{q}=1$ and $p,z,q> 1$. We treat each of the previous terms separately:
\begin{enumerate}
\item[$\bullet$] First observe that for $1<p<n/(1-\omega)$ we have for the term $(1)$ above:
\begin{equation*}
\||x-x(s_0)|^{\omega-1}\|_{L^p(B_r)}\leq C r^{n/p+\omega-1}.
\end{equation*}
\item[$\bullet$] By hypothesis $v(\cdot, s_0)\in bmo$ and applying the Lemma \ref{PropoBMO1} we have $\|v-\overline{v}_{B_r}\|_{L^z(B_r)}\leq C|B_r|^{1/z}\|v(\cdot,s_0)\|_{bmo}$. Since $\underset{s_0\in [0,T]}{\sup}\; \|v(\cdot,s_0)\|_{bmo}\leq \mu$, we write for the term $(2)$:
\begin{equation*}
\|v-\overline{v}_{B_r}\|_{L^z(B_r)}\leq C\mu\; r^{n/z}.
\end{equation*}
\item[$\bullet$] Finally for $(3)$ by the maximum principle we have $\|\psi_+(\cdot, s_0)\|_{L^q(B_r)}\leq  \|\psi_+(\cdot, 0)\|_{L^q}$; hence using the fact that $\psi_0$ is an $r$-molecule and remark \ref{Remark2}-2) we obtain
\begin{equation*}
\|\psi_+(\cdot, s_0)\|_{L^q(B_r)}\leq  C \bigg[r^{-\gamma}\bigg]^{1/q}\left[\frac{1}{r^{n+\gamma}}\right]^{1-1/q}.
\end{equation*}
\end{enumerate}
We combine all these inequalities together in order to obtain the following estimation for (\ref{Equa1}):
\begin{equation}\label{Bola1}
I_{1,B_r}\leq C\mu\; r^{\omega-1-\gamma}.
\end{equation}
\item[(ii)] \underline{Estimations for the dyadic corona $E_k$}. Let us note $I_{1,E_k}$ the integral 
$$I_{1,E_k}=\int_{E_k}|x-x(s_0)|^{\omega-1}|v-\overline{v}_{B_r}| |\psi_+(x,s_0)|dx.$$
Since over $E_k$ we have\footnote{recall that $0<\gamma<\omega<2\delta\leq 1$.} $|x-x(s_0)|^{\omega-1}\leq C 2^{k(\omega-1)}r^{\omega-1}$ we write
\begin{eqnarray*}
I_{1,E_k}&\leq & C2^{k(\omega-1)}r^{\omega-1}\left(\int_{E_k}|v-\overline{v}_{B_{r2^k}}| |\psi_+(x,s_0)|dx+\int_{E_k}|\overline{v}_{B_r}-\overline{v}_{B_{r2^k}}| |\psi_+(x,s_0)|dx\right)
\end{eqnarray*}
where we noted $B_{r2^k}=B(x(s_0),r2^k)$, then
\begin{eqnarray*}
I_{1,E_k} &\leq& C2^{k(\omega-1)}r^{\omega-1}\left(\int_{B_{r2^k}}|v-\overline{v}_{B_{r2^k}}| |\psi_+(x,s_0)|dx+\int_{B_{r2^k}}|\overline{v}_{B_r}-\overline{v}_{B_{r2^k}}| |\psi_+(x,s_0)|dx\right).
\end{eqnarray*}
Now, since $v(\cdot, s_0)\in bmo(\mathbb{R}^n)$, using the Lemma \ref{PropoBMO1} we have $|\overline{v}_{B_r}-\overline{v}_{B_{r2^k}}| \leq Ck\|v(\cdot,s_0)\|_{bmo}\leq Ck\mu$ and we can write
\begin{eqnarray*}
I_{1,E_k}&\leq &C2^{k(\omega-1)}r^{\omega-1}\left(\int_{B_{r2^k}}|v-\overline{v}_{B_{r2^k}}| |\psi_+(x,s_0)|dx+Ck\mu\|\psi_+(\cdot,s_0)\|_{L^1}\right)\\[5mm]
&\leq & C2^{k(\omega-1)}r^{\omega-1}\left(\|\psi_+(\cdot,s_0)\|_{L^{a_0}} \|v-\overline{v}_{B_{r2^k}}\|_{L^{\frac{a_0}{a_0-1}}} +Ck\mu\;  r^{-\gamma}\right)
\end{eqnarray*}
where we used Hölder's inequality with $1<a_0<\frac{n}{n+(\omega-1)}$ and maximum principle for the last term above. Using again the properties of $bmo$ spaces we have
$$I_{1,E_k}\leq  C2^{k(\omega-1)}r^{\omega-1}\left(\|\psi_+(\cdot,0)\|_{L^1}^{1/a_0}\|\psi_+(\cdot,0)\|_{L^\infty}^{1-1/a_0} |B_{r2^k}|^{1-1/a_0}\|v(\cdot,s)\|_{bmo} +Ck\mu r^{-\gamma}\right).$$
Let us now apply the estimates given by hypothesis for $\|\psi_+(\cdot,0)\|_{L^1}$, $\|\psi_+(\cdot,0)\|_{L^\infty}$ and $\|v(\cdot,s_0)\|_{bmo}$  to obtain
$$I_{1,E_k}\leq  C2^{k(n-n/a_0+\omega-1)}r^{\omega-1-\gamma}\mu +C2^{k(\omega-1)}k\mu\;  r^{\omega-1-\gamma}.$$
Since $1<a_0<\frac{n}{n+(\omega-1)}$, we have $n-n/a_0+(\omega-1)<0$, so that, summing over each dyadic corona $E_k$, we have
\begin{equation}\label{Coronak}
\sum_{k\geq 1}I_{1,E_k}\leq C\mu\; r^{\omega-1-\gamma}.
\end{equation}
\end{enumerate}
Finally, gathering together the estimations (\ref{Bola1}) and (\ref{Coronak}) we obtain the desired conclusion.\hfill$\blacksquare$\\

\begin{Lemme}\label{Lemme2}
For integral $I_2$ in inequality (\ref{smallEstrella}) we have the estimate $I_2\leq C r^{\omega-1-\gamma}$.
\end{Lemme}
\textbf{\textit{Proof.}} As for the Lemma \ref{Lemme1}, we consider $\mathbb{R}^n$ as the union of a ball with dyadic coronas centered on $x(s_0)$ (cf. (\ref{SmallDecoupage})). 
\begin{enumerate}
\item[(i)] \underline{Estimations over the ball $B_r$}. We write, using the maximum principle and the hypothesis for $\|\psi_+(\cdot, 0)\|_{L^\infty}$:
\begin{eqnarray*}
I_{2,B_r}&=&\int_{B_r}\big|\mathcal{L}(|x-x(s_0)|^{\omega})\big||\psi_+(x,s_0)|dx \leq \|\psi_+(\cdot, s_0)\|_{L^{\infty}}\int_{B_r}|\mathcal{L}|x-x(s_0)|^{\omega}|dx\nonumber\\
&\leq & \|\psi_+(\cdot, 0)\|_{L^{\infty}}\int_{\{|x|\leq r\}}\left|\mbox{v.p.}\int_{\mathbb{R}^n}[|x|^{\omega}-|x-y|^{\omega}]\pi(y)dy\right|dx\nonumber\\
&\leq & r^{-n-\gamma}\int_{\{|x|\leq r\}}\left|\mbox{v.p.}\int_{\mathbb{R}^n}[|x|^{\omega}-|x-y|^{\omega}]\pi(y)dy\right|dx. 
\end{eqnarray*}
We use now the hypothesis (\ref{DefKernel2}) and (\ref{DefKernel3}) for the function $\pi$ in the case \textbf{(c)}, \textit{i.e.} $\alpha=\beta=1/2$ and  $0<\delta<1/2$, in order to obtain
\begin{eqnarray}
I_{2,B_r}&\leq&  r^{-n-\gamma}\int_{\{|x|\leq r\}}\left|\mbox{v.p.}\int_{\{|y|\leq 1\}}\frac{|x|^{\omega}-|x-y|^{\omega}}{|y|^{n+1}}dy\right|dx+ r^{-n-\gamma}\int_{\{|x|\leq r\}}\int_{\mathbb{R}^n}\frac{\left||x|^{\omega}-|x-y|^{\omega}\right|}{|y|^{n+2\delta}}dydx\nonumber\\
& \leq & r^{-n-\gamma}\big( I_{2,B_r^1}+I_{2,B_r^{\delta}}\big). \label{Estimate101}
\end{eqnarray}
We start studying the first term $I_{2,B_r^1}$ above. Recalling that 
\begin{equation}\label{DefDeltaOmega}
(-\Delta)^{1/2}(|x|^\omega)=\mbox{v.p.}\displaystyle{\int_{\mathbb{R}^n}}\frac{|x|^{\omega}-|x-y|^{\omega}}{|y|^{n+1}}dy=|x|^{\omega-1},
\end{equation}
 by homogeneity and using the fact that $0<r<1$ we obtain:
$$ I_{2,B_r^1}\leq r^{\omega+n-1}\left(\int_{\{|x|\leq 1\}}|x|^{\omega-1}dx+\int_{\{|x|\leq 1\}}\int_{\{|y|>1/r\}}\frac{||x|^{\omega}-|x-y|^{\omega}|}{|y|^{n+1}}dydx\right)=Cr^{\omega+n-1}.$$
For the second term $I_{2,B_r^{\delta}}$ we will proceed as follows. First, by homogeneity we obtain
\begin{equation*}
I_{2,B^\delta_r}=r^{\omega+n-2\delta}\underbrace{\int_{\{|x|\leq 1\}}\int_{\mathbb{R}^n}\frac{\left||x|^{\omega}-|x-y|^{\omega}\right|}{|y|^{n+2\delta}}dydx}_{I}.
\end{equation*}
Then we decompose this integral $I$ in the following way
\begin{eqnarray*}
I&=&\int_{\{|x|\leq 1\}}\int_{\{|y|\leq 1\}}\frac{\left||x|^{\omega}-|x-y|^{\omega}\right|}{|y|^{n+2\delta}}dydx+\int_{\{|x|\leq 1\}}\int_{\{|y|>1\}}\frac{\left||x|^{\omega}-|x-y|^{\omega}\right|}{|y|^{n+2\delta}}dydx\\
&\leq & \int_{\{|x|\leq 1\}}\left(\underset{0<|y|<1}{\sup}\frac{\left||x|^{\omega}-|x-y|^{\omega}\right|}{|y|}\right) \left(\int_{\{|y|\leq 1\}}|y|^{1-n-2\delta}dy\right)dx+\int_{\{|x|\leq 1\}}\left(\int_{\{|y|>1\}} |y|^{\omega-n-2\delta}dy\right)dx
\end{eqnarray*}
Since $0<\gamma<\omega<2\delta< 1$, it is not complicated to see that 
\begin{eqnarray}\label{FirstHomogeneityEstimate}
I &\leq & C\int_{\{|x|\leq 1\}}\left(\underset{0<|y|<1}{\sup}\frac{\left||x|^{\omega}-|x-y|^{\omega}\right|}{|y|}\right)dx+ C
\end{eqnarray}
and that this latter quantity is bounded. Then, getting back to (\ref{Estimate101}) we write $I_{2,B_r}\leq C (r^{\omega-\gamma-1}+r^{\omega-\gamma-2\beta}).$ Recalling that we are working with small molecules, \textit{i.e.} that $0<r<1$ , we obtain $r^{\omega-2\beta-\gamma}\leq r^{\omega-1-\gamma}$, so we finally have
\begin{equation*}
I_{2,B_r}\leq C r^{\omega-\gamma-1}.
\end{equation*}
The case \textbf{(d)}, when  $\alpha=\beta=\delta=1/2$, is easier since $(-\Delta)^{1/2}(|x|^\omega)=|x|^{\omega-1}$. Thus, in any case we can write:
\begin{equation}\label{smallBola2}
I_{2,B_r}=\int_{B_r}|\mathcal{L}(|x-x(s_0)|^{\omega})||\psi_+(x,s_0)|dx \leq Cr^{\omega-1-\gamma}.
\end{equation}
\item[(ii)]  \underline{Estimations for the dyadic corona $E_k$}. We start with the case \textbf{(c)} when $\alpha=\beta=1/2$ and  $0<\delta<1/2$:
\begin{eqnarray*}
I_{2, E_k}&=&\int_{E_k}|\mathcal{L}(|x-x(s_0)|^{\omega})| |\psi_+(x,s_0)|dx  \leq \|\psi_+(\cdot, s_0)\|_{L^1}\underset{x\in E_k}{\sup}\left|\mathcal{L}(|x-x(s_0)|^{\omega})\right| \\
&\leq & r^{-\gamma} \bigg(\underbrace{\underset{r2^{k-1}<|x|\leq r2^{k}}{\sup}\left|\mbox{v.p.}\int_{\{|y|\leq 1\}}\frac{|x|^{\omega}-|x-y|^{\omega}}{|y|^{n+1}}dy\right|}_{I_{2,E_k^1}} +\underbrace{\underset{r2^{k-1}<|x|\leq r2^{k}}{\sup}\int_{\mathbb{R}^n}\frac{||x|^{\omega}-|x-y|^{\omega}|}{|y|^{n+2\delta}}dy}_{I_{2,E_k^\delta}} \bigg)
\end{eqnarray*}
Let us start with $I_{2,E_k^1}$, by homogeneity and using the formula (\ref{DefDeltaOmega}) we obtain
\begin{eqnarray*}
I_{2,E_k^1}&\leq& \underset{r2^{k-1}<|x|\leq r2^{k}}{\sup}|x|^{\omega-1} + C(r2^{k})^{\omega-1}\bigg(\underset{1<|x|\leq 2}{\sup}\int_{\{|y|>1/r2^{k-1}\}}\frac{||x|^{\omega}-|x-y|^{\omega}|}{|y|^{n+1}}dy\bigg)
\end{eqnarray*}
We only need to study the last term of this expression. If $0<r2^{k-1}\leq 1$, the integral above is immediately bounded by a constant. The case when $1<r2^{k-1}$ is treated as follows:
\begin{eqnarray*}
\underset{1<|x|\leq 2}{\sup}\int_{\{|y|>1/r2^{k-1}\}}\frac{||x|^{\omega}-|x-y|^{\omega}|}{|y|^{n+1}}dy&=&\underset{1<|x|\leq 2}{\sup}\bigg(\int_{\{1/r2^{k-1}<|y|<1\}}\frac{||x|^{\omega}-|x-y|^{\omega}|}{|y|^{n+1}}dy+\underset{\{1<|y|\}}{\int}\frac{||x|^{\omega}-|x-y|^{\omega}|}{|y|^{n+1}}dy\bigg)\\
&\leq & \underset{1<|x|\leq 2}{\sup}\bigg(\underset{0<|y|<1}{\sup}\frac{||x|^{\omega}-|x-y|^{\omega}|}{|y|}\bigg)\ln(2^{k-1})+C
\end{eqnarray*}
Thus we obtain $I_{2,E_k^1}\leq C(r2^{k})^{\omega-1}\big(1+\ln(2^{k-1})\big)$.\\

The term $I_{2,E_k^\delta}$ is easier: applying essentially the same ideas used in the formulas (\ref{Estimate101})-(\ref{FirstHomogeneityEstimate}) above and by homogeneity we have $I_{2,E_k^\delta}\leq C (r2^k)^{\omega-2\delta}$.\\

Finally, we obtain the following inequality for $I_{2,E_k}$:
\begin{eqnarray*}
I_{2,E_k}& \leq &  Cr^{-\gamma}\left( (r2^k)^{\omega-1}\big(1+\ln(2^{k-1})\big)+(r2^k)^{\omega-2\delta}\right)
\end{eqnarray*}
Since $0<\gamma<\omega<2\delta<1$, summing over $k\geq 1$, we obtain $\displaystyle{\sum_{k\geq 1}}I_{2,E_k} \leq  Cr^{-\gamma}\left( r^{\omega-1}+r^{\omega-2\delta}\right)$. Repeating the same argument used before (\textit{i.e.} the fact that $0<r<1$), we finally obtain
\begin{equation}\label{smallCoronak2}
\sum_{k\geq 1}I_{2,E_k}\leq C r^{\omega-1-\gamma}.
\end{equation}
\end{enumerate}
The case \textbf{(d)} is straightforward since we have $\mathcal{L}=(-\Delta)^{1/2}$ and $(-\Delta)^{1/2}(|x|^\omega)=|x|^{\omega-1}$.\\

In order to finish the proof of Lemma \ref{Lemme2} we combine together the estimates (\ref{smallBola2}) and (\ref{smallCoronak2}).\hfill$\blacksquare$\\

Now we continue the proof of the Proposition \ref{Propo1}. Using the Lemmas \ref{Lemme1} and \ref{Lemme2} and getting back to estimate (\ref{smallEstrella}) we have
$$\left|\partial_{s_0} \int_{\mathbb{R}^n}\Omega_0(x)\psi_+(x,s_0)dx\right| \leq  C(\mu+1)\; r^{\omega-1-\gamma}$$

This last estimation is compatible with the estimate (\ref{SmallConcentration}) for $0\leq s_0\leq \epsilon r$ small enough: just fix $K$ such that
\begin{equation}\label{SmallConstants}
C\left(\mu+1\right)\leq K(\omega-\gamma).
\end{equation}
Indeed, since the time $s_0$ is very small, we can linearize the formula $(r+Ks_0)^{\omega-\gamma}$ in the right-hand side of (\ref{SmallConcentration}) in order to obtain
\begin{equation*}
\phi=(r+Ks_0)^{\omega-\gamma} \thickapprox r^{\omega-\gamma}\left(1+[K(\omega-\gamma)]\frac{s_0}{r}\right).
\end{equation*}
Finally, taking the derivative with respect to $s_0$ in the above expression we have $\phi' \thickapprox r^{\omega-1-\gamma}K(\omega-\gamma)$ and with condition (\ref{SmallConstants}) Proposition \ref{Propo1} follows.\hfill$\blacksquare$\\

Now we will give a sligthly different proof of the maximum principle of A. C\'ordoba \& D. C\'ordoba. Indeed, the following proof only relies on the concentration condition proved in the lines above. 
\begin{Proposition}[Small time Height condition]\label{Propo2}
Under the  hypothesis of Theorem \ref{SmallGeneralisacion}, if $\psi(x,s_0)$ satisfies the concentration condition (\ref{SmallConcentration}), then we have the following height condition
\begin{equation}\label{SmallForCondHauteur}
\|\psi(\cdot, s_0)\|_{L^\infty}  \leq \frac{1}{\left(r+Ks_0\right)^{n+\gamma}}.
\end{equation}
\end{Proposition}
\textit{\textbf{Proof.}} Assume that molecules we are working with are smooth enough. Following an idea of \cite{Cordoba} (section 4 p.522-523) (see also \cite{Jacob} p. 346), we will note $\overline{x}$ the point of $\mathbb{R}^n$ such that $\psi(\overline{x},s_0)=\|\psi(\cdot,s_0)\|_{L^\infty}$. Thus we can write, by the properties of the function $\pi$ (recall that we assumed $\alpha=\beta=1/2$ and $0<\delta<1/2$ or $\alpha=\beta=\delta=1/2$):
\begin{equation}\label{Infty1}
\frac{d}{ds_0}\|\psi(\cdot,s_0)\|_{L^\infty}\leq -\int_{\mathbb{R}^n}|\psi(\overline{x},s_0)-\psi(\overline{x}-y,s_0)|\pi(y)dy\leq -\int_{\{|\overline{x}-y|<1\}}\frac{|\psi(\overline{x},s_0)-\psi(y,s_0)|}{|\overline{x}-y|^{n+1}}dy\leq 0.
\end{equation}
For simplicity, we will assume that $\psi(\overline{x},s_0)$ is positive. Let us consider the corona centered in $\overline{x}$ defined by
$$\mathcal{C}(R_1,R_2)=\{y\in \mathbb{R}^n:R_1\leq|\overline{x}-y|\leq R_2\}$$ 
where $1>R_2=\rho R_1$ with $\rho >2$ and where $R_1$ will be fixed later. Then:
\begin{equation*}
\int_{\{|\overline{x}-y|<1\}}\frac{\psi(\overline{x},s_0)-\psi(y,s_0)}{|\overline{x}-y|^{n+1}}dy\geq \int_{\mathcal{C}(R_1,R_2)}\frac{\psi(\overline{x},s_0)-\psi(y,s_0)}{|\overline{x}-y|^{n+1}}dy.
\end{equation*}
Define the sets $B_1$ and $B_2$ by $B_1=\{y\in \mathcal{C}(R_1,R_2): \psi(\overline{x},s_0)-\psi(y,s_0)\geq \frac{1}{2}\psi(\overline{x},s_0)\}$ and $B_2=\{y\in \mathcal{C}(R_1,R_2): \psi(\overline{x},s_0)-\psi(y,s_0)< \frac{1}{2}\psi(\overline{x},s_0)\}$ such that $\mathcal{C}(R_1,R_2)=B_1\cup B_2$. \\

We obtain the inequalities
\begin{eqnarray*}
\int_{\mathcal{C}(R_1,R_2)}\frac{\psi(\overline{x},s_0)-\psi(y,s_0)}{|\overline{x}-y|^{n+1}}dy &\geq &\int_{B_1}\frac{\psi(\overline{x},s_0)-\psi(y,s_0)}{|\overline{x}-y|^{n+1}}dy \geq \frac{\psi(\overline{x},s_0)}{2R_2^{n+1}}|B_1|=\frac{\psi(\overline{x},s_0)}{2R_2^{n+1}}\left(|\mathcal{C}(R_1,R_2)|-|B_2|\right).
\end{eqnarray*}
Since $R_2=\rho R_1$ one has
\begin{equation}\label{Infty3}
\int_{\mathcal{C}(R_1,R_2)}\frac{\psi(\overline{x},s_0)-\psi(y,s_0)}{|\overline{x}-y|^{n+1}}dy\geq \frac{\psi(\overline{x},s_0)}{2\rho^{n+1}R_1^{n+1}}\bigg(v_n(\rho^n -1)R_1^n-|B_2|\bigg)
\end{equation}
where $v_n$ denotes the volume of the $n$-dimensional unit ball.\\ 

Now, we will estimate the quantity $|B_2|$ in terms of $\psi(\overline{x},s_0)$ and $R_1$ with the following lemma.
\begin{Lemme}
For the set $B_2$ we have the following estimations
\begin{enumerate}
\item[1)] if $|\overline{x}-x(s_0)|>2R_2$ then $C_1(r+Ks_0)^{\omega-\gamma}\psi(\overline{x},s_0)^{-1}R_1^{-\omega}\geq |B_2|$.\\

\item[2)] if $|\overline{x}-x(s_0)|<R_1/2$ then $C_1(r+Ks_0)^{\omega-\gamma}\psi(\overline{x},s_0)^{-1}R_1^{-\omega}\geq |B_2|$.\\

\item[3)] if $R_1/2\leq |\overline{x}-x(s_0)|\leq 2R_2$ then $\big(C_2 (r+Ks_0)^{\omega-\gamma}R_1^{n-\omega}\psi(\overline{x},s_0)^{-1}\big)^{1/2}\geq |B_2|$.
\end{enumerate}
\end{Lemme}
Recall that for the molecule's center $x_0\in \mathbb{R}^n$ we noted its transport by $x(s_0)$ which is defined by formula (\ref{Defpointx_0}).\\

\textit{\textbf{Proof.}} For all these estimates, our starting point is the concentration condition (\ref{SmallConcentration}):
\begin{eqnarray}
(r+Ks_0)^{\omega-\gamma}\geq \int_{\mathbb{R}^n}|\psi(y,s_0)||y-x(s_0)|^{\omega}dy &\geq & \int_{B_2}|\psi(y,s_0)||y-x(s_0)|^{\omega}dy \nonumber\geq \frac{\psi(\overline{x},s_0)}{2}\int_{B_2}|y-x(s_0)|^{\omega}dy.\label{EstimationB2}
\end{eqnarray}
We just need to estimate the last integral following the cases given by the lemma. The first two cases are very similar. Indeed, if $|\overline{x}-x(s_0)|>2R_2$  then we have
$$\underset{y\in B_2\subset \mathcal{C}(R_1,R_2)}{\min}|y-x(s_0)|^{\omega}\geq R_2^{\omega}=\rho^{\omega}R_1^{\omega}$$
while for the second case, if $|\overline{x}-x(s_0)|<R_1/2$, one has
$$\underset{y\in B_2\subset \mathcal{C}(R_1,R_2)}{\min}|y-x(s_0)|^{\omega}\geq \frac{R_1^{\omega}}{2^\omega}.$$
Applying these results to (\ref{EstimationB2}) we obtain $(r+Ks_0)^{\omega-\gamma}\geq \frac{\psi(\overline{x},s_0)}{2} \rho^{\omega} R_1^{\omega}|B_2|$ and  $(r+Ks_0)^{\omega-\gamma}\geq \frac{\psi(\overline{x},s_0)}{2} \frac{R_1^{\omega}}{2^\omega}|B_2|$, and since $\rho>2$ we have the desired estimate
\begin{equation}\label{FormEstima1}
\frac{C_1 (r+Ks_0)^{\omega-\gamma}}{\psi(\overline{x},s_0) R_1^{\omega}} \geq \frac{2 (r+Ks_0)^{\omega-\gamma}}{\rho^{\omega}\psi(\overline{x},s_0) R_1^{\omega}} \geq |B_2| \qquad \mbox{with } C_1=2^{1+\omega}.
\end{equation}
For the last case, since $R_1/2\leq |\overline{x}-x(s_0)|\leq 2R_2$ we can write using the Cauchy-Schwarz inequality
\begin{equation}\label{HolderInver}
\int_{B_2}|y-x(s_0)|^{\omega}dy\geq |B_2|^2\left(\int_{B_2}|y-x(s_0)|^{-\omega}dy\right)^{-1}
\end{equation}
Now, observe that in this case we have $B_2\subset B(x(s_0), 5R_2)$ and then
$$\int_{B_2}|y-x(s_0)|^{-\omega}dy\leq \int_{B(x(s_0), 5 R_2)}|y-x(s_0)|^{-\omega}dy\leq v_n (5\rho R_1)^{n-\omega}.$$
Getting back to(\ref{HolderInver}) we obtain
$$\int_{B_2}|y-x(s_0)|^{\omega}dy\geq |B_2|^2 v_n^{-1} (5 \rho R_1)^{-n+\omega}$$
We use this estimate in (\ref{EstimationB2}) to obtain
\begin{equation}\label{FormEstima2}
\frac{C_2(r+Ks_0)^{\omega/2-\gamma/2} R_1^{n/2-\omega/2}}{\psi(\overline{x},s_0)^{1/2}}\geq |B_2|,
\end{equation}
where $C_2=(2\times 5^{n-\omega} v_n\rho^{n-\omega})^{1/2}$. The lemma is proven.\hfill $\blacksquare$\\

With estimates (\ref{FormEstima1}) and (\ref{FormEstima2}) at our disposal we can write
\begin{enumerate}
\item[(i)] if $|\overline{x}-x(s_0)|>2R_2$ or $|\overline{x}-x(s_0)|<R_1/2$ then
\begin{equation*}
\int_{\mathcal{C}(R_1,R_2)}\frac{\psi(\overline{x},s_0)-\psi(y,s_0)}{|\overline{x}-y|^{n+1}}dy\geq  \frac{\psi(\overline{x},s_0)}{2\rho^{n+1}R_1^{n+1}}\bigg(v_n(\rho^n -1)R_1^n-\frac{C_1(r+Ks_0)^{\omega-\gamma}}{\psi(\overline{x},s_0)} R_1^{-\omega}\bigg)
\end{equation*}
\item[(ii)] if $R_1/2\leq |\overline{x}-x(s_0)|\leq 2R_2$
\begin{equation*}
\int_{\mathcal{C}(R_1,R_2)}\frac{\psi(\overline{x},s_0)-\psi(y,s_0)}{|\overline{x}-y|^{n+1}}dy\geq  \frac{\psi(\overline{x},s_0)}{2\rho^{n+1}R_1^{n+1}}\bigg(v_n(\rho^n -1)R_1^n-\frac{C_2 (r+Ks_0)^{\omega/2-\gamma/2}R_1^{n/2-\omega/2}}{\psi(\overline{x},s_0)^{1/2}}\bigg)
\end{equation*}
\end{enumerate}
Now, if we set $R_1=(r+Ks_0)^{\frac{(\omega-\gamma)}{n+\omega}}\psi(\overline{x},s_0)^{\frac{-1}{n+\omega}}$ and if $\rho$ is big enough such that the expression in brackets above is positive, we obtain for cases (i) and (ii) the following estimate for (\ref{Infty3}):
\begin{equation*}
\int_{\mathcal{C}(R_1,R_2)}\frac{\psi(\overline{x},s_0)-\psi(y,s_0)}{|\overline{x}-y|^{n+1}}dy\geq C (r+Ks_0)^{-\frac{(\omega-\gamma)}{n+\omega}} \psi(\overline{x},s_0)^{1+\frac{1}{n+\omega}}
\end{equation*}
where $C=C(n,\rho)=\frac{v_n (\rho^n-1)-\sqrt{2v_n}(5\rho)^{\frac{n-\omega}{2}}}{2\rho^{n+1}}<1$ is a small positive constant. Now, and for all possible cases considered before, we have the following estimate for (\ref{Infty1}):
$$\frac{d}{ds_0}\|\psi(\cdot,s_0)\|_{L^\infty}\leq -C(r+Ks_0)^{-\frac{(\omega-\gamma)}{n+\omega}} \|\psi(\cdot,s_0)\|_{L^\infty}^{1+\frac{1}{n+\omega}}.$$
Solving this differential inequality with initial data $\|\psi(\cdot, 0)\|_{L^{\infty}}\leq r^{-n-\gamma}$, we obtain $\|\psi(\cdot,s_0)\|_{L^\infty}\leq(r+Ks_0)^{-(n+\gamma)}$.\\ 

The proof of Proposition \ref{Propo2} is finished for regular molecules. In order to obtain the global result, remark that, for viscosity solutions (\ref{SistApprox}), we have that $\Delta \theta(\overline{x},s_0)\leq 0$ at the points $\overline{x}$ where $\theta(\cdot,s_0)$ reaches its maximum value. See \cite{Cordoba} for more details. \hfill$\blacksquare$\\

We treat now the last part of Theorem \ref{SmallGeneralisacion}:
\begin{Proposition}[First $L^1$ estimate]\label{Propo3}
If $\psi(x,s_0)$ is a solution of the problem (\ref{SmallEvolution}), then we have the following $L^1$-norm estimate:
$$\|\psi(\cdot, s_0)\|_{L^1}\leq  \frac{v_n }{\big(r+Ks_0\big)^{\gamma}}.$$
\end{Proposition}
\textit{\textbf{Proof.}} We write
\begin{eqnarray}
\int_{\mathbb{R}^n}|\psi(x,s_0)|dx&=&\int_{\{|x-x(s_0)|< D\}}|\psi(x,s_0)|dx+\int_{\{|x-x(s_0)|\geq D\}}|\psi(x,s_0)|dx\nonumber\\
&\leq & v_n D^n \|\psi(\cdot, s_0)\|_{L^\infty}+D^{-\omega}\int_{\mathbb{R}}|\psi(x,s_0)||x-x(s_0)|^\omega dx\nonumber
\end{eqnarray}
Now using (\ref{SmallForCondHauteur}) and (\ref{SmallConcentration}) one has:
\begin{eqnarray*}
\int_{\mathbb{R}^n}|\psi(x,s_0)|dx&\leq & v_n \frac{D^n }{\left(r+Ks_0\right)^{n+\omega}}  +D^{-\omega}(r+Ks_0)^{\omega-\gamma}
\end{eqnarray*}
where $v_n$ denotes the volume of the unit ball. To continue, it is enough to choose correctly the real parameter $D$ to obtain
$$\|\psi(\cdot, s_0)\|_{L^1}\leq \frac{v_n}{\big(r+Ks_0\big)^{\gamma}}$$
\begin{flushright}$\blacksquare$\end{flushright}
\subsection{Molecule's evolution: Second step}\label{SecEvolMol2}
In the previous section we have obtained deformed molecules after a very small time $s_0$. The next theorem shows us how to obtain similar profiles in the inputs and the outputs in order to perform an iteration in time. \\

Recall that we consider here a Lévy-type operator $\mathcal{L}$ of the form (\ref{DefOperator1}) with an associate Lévy measure $\pi$ that satisfies hypothesis (\ref{DefKernel2}) and (\ref{DefKernel3}) with the following values of the parameters $\alpha, \beta, \delta$:
\begin{itemize}
\item[\textbf{(c)}] $\alpha=\beta=1/2$ and $0<\delta<1/2$,
\item[\textbf{(d)}] $\alpha=\beta=\delta=1/2$.
\end{itemize}
\begin{Theoreme}\label{Generalisacion} Set $\gamma$ and $\omega$ two real numbers such that $0<\gamma<\omega<2\delta<1$ in the case \textbf{(c)} or $0<\gamma<\omega<1$ in the case \textbf{(d)}. Let $0< s_1\leq T$ and let $\psi(x,s_1)$ be a solution of the problem
\begin{equation}\label{Evolution}
\left\lbrace
\begin{array}{rl}
\partial_{s_1} \psi(x,s_1)=& -\nabla\cdot(v\, \psi)(x,s_1)-\mathcal{L}\psi(x,s_1)\\[5mm]
\psi(x,0)=& \psi(x,s_0)  \qquad \qquad \mbox{with } s_0>0\\[5mm]
div(v)=&0 \quad \mbox{and }\; v\in L^{\infty}([0,T];bmo(\mathbb{R}^n))\quad \mbox{with } \underset{s_1\in [s_0,T]}{\sup}\; \|v(\cdot,s_1)\|_{bmo}\leq \mu
\end{array}
\right.
\end{equation}
If $\psi(x,s_0)$ satisfies the three following conditions
\begin{eqnarray*}
\int_{\mathbb{R}^n}|\psi(x,s_0)||x-x(s_0)|^\omega dx \leq (r+Ks_0)^{\omega-\gamma}; \quad \|\psi(\cdot, s_0)\|_{L^\infty}\leq \frac{1}{\left(r+ Ks_0\right)^{n+\gamma }}; \quad \|\psi(\cdot, s_0)\|_{L^1} \leq  \frac{v_n}{\big(r+Ks_0\big)^{\gamma }}
\end{eqnarray*}
where $K=K(\mu)$ is given by (\ref{SmallConstants}) and $s_0$ is such that $(r+Ks_0)<1$. Then for all $0< s_1\leq\epsilon r$ small, we have the following estimates
\begin{eqnarray}
\int_{\mathbb{R}^n}|\psi(x,s_1)||x-x(s_1)|^\omega dx &\leq &(r+K(s_0+s_1))^{\omega-\gamma}  \label{Concentration2}\\
\|\psi(\cdot,s_1)\|_{L^\infty}&\leq & \frac{1}{\left(r+K(s_0+s_1)\right)^{n+\gamma }}\label{Linftyevolutionnotsmalltime}\\
\|\psi(\cdot,s_1)\|_{L^1} &\leq & \frac{v_n}{\big(r+K(s_0+s_1)\big)^{\gamma }} \label{L1evolutionsmalltime}
\end{eqnarray}
\end{Theoreme}

\begin{Remarque}
\begin{itemize}
\item[]
\item[1)] Since $s_1$ is small and $(r+Ks_0)<1$, we can without loss of generality assume that $(r+K(s_0+s_1))<1$: otherwise, by the maximum principle there is nothing to prove.
\item[2)] The new molecule's center $x(s_1)$ used in formula (\ref{Concentration2}) is fixed by 
\begin{equation}\label{Defpointx_s}
\left\lbrace
\begin{array}{rl}
x'(s_1)=& \overline{v}_{B_{f_1}}=\frac{1}{|B_{f_1}|}\displaystyle{\int_{B_{f_1}}}v(y,s_1)dy\\[5mm]
x(0)=& x(s_0).
\end{array}
\right.
\end{equation}  
And here we noted $B_{f_1}=B(x(s_1),f_1)$ with $f_1$ a real valued function given by 
\begin{equation}\label{DefiFunctionF}
f_1=(r+Ks_0).
\end{equation}
Note that by remark $1)$ above we have $0<f_1<1$.
\end{itemize}
\end{Remarque}
We will follow the same scheme as before:  we prove the concentration condition (\ref{Concentration2}), with this estimate at hand we will control the $L^\infty$ decay in Proposition \ref{Propo5} and then we will obtain the suitable $L^1$ control in Proposition \ref{Propo6}.
\begin{Proposition}[Concentration condition]\label{Propo4}
Under the  hypothesis of Theorem \ref{Generalisacion}, if $\psi(\cdot,s_0)$ is an initial data then the solution $\psi(x,s_1)$ of (\ref{Evolution}) satisfies
\begin{equation*}
\int_{\mathbb{R}^n} |\psi(x,s_1)||x-x(s_1)|^{\omega}dx \leq (r+K(s_0+s_1))^{\omega-\gamma}
\end{equation*}
for $x(s_1)\in \mathbb{R}^n$ given by formula (\ref{Defpointx_s}), with $0\leq s_1\leq \epsilon r$. 
\end{Proposition}
\textit{\textbf{Proof}.} The calculations are very similar of those of Proposition \ref{Propo1}: the only difference stems from the initial data and the definition of the center $x(s_1)$. So, let us write $\Omega_1(x)=|x-x(s_1)|^{\omega}$ and $\psi(x)=\psi_+(x)-\psi_-(x)$ where the functions $\psi_{\pm}(x)\geq 0$ have disjoint support. 
Thus, by linearity and using the positivity principle we have  
$$|\psi(x,s_1)|=|\psi_+(x,s_1)-\psi_-(x,s_1)|\leq \psi_+(x,s_1)+\psi_-(x,s_1)$$ and we can write
$$\int_{\mathbb{R}^n}|\psi(x,s_1)\Omega_1(x)dx\leq\int_{\mathbb{R}^n}\psi_+(x,s_1)\Omega_1(x)dx+\int_{\mathbb{R}^n}\psi_-(x,s_1)\Omega_1(x)dx$$
so we only have to treat one of the integrals on the right-hand side above. We have:
\begin{eqnarray*}
I&=&\left|\partial_{s_1} \int_{\mathbb{R}^n}\Omega_1(x)\psi_+(x,s_1)dx\right|=\left|\int_{\mathbb{R}^n}-\nabla\Omega_1(x)\cdot x'(s_1)\psi_+(x,s_1)+\Omega_1(x)\left[-\nabla\cdot(v\, \psi_+(x,s_1))-\mathcal{L}\psi_+(x,s_1)\right]dx\right|
\end{eqnarray*}
Using the fact that $v$ is divergence free, we obtain
\begin{equation*}
I=\left|\int_{\mathbb{R}^n}\nabla\Omega_1(x)\cdot(v-x'(s_1))\psi_+(x,s_1)-\Omega_1(x)\mathcal{L}\psi_+(x,s_1)dx\right|.
\end{equation*}
Finally, using the definition of $x'(s_1)$ given in (\ref{Defpointx_s}) and replacing $\Omega_1(x)$ by $|x-x(s_1)|^{\omega}$ in the first integral we obtain
\begin{equation}\label{Estrella}
I\leq c \underbrace{\int_{\mathbb{R}^n}|x-x(s_1)|^{\omega-1}|v-\overline{v}_{B_{f_1}}| |\psi_+(x,s_1)|dx}_{I_1} + c\underbrace{\int_{\mathbb{R}^n}|\mathcal{L}\Omega_1(x)||\psi_+(x,s_1)|dx}_{I_2}.
\end{equation}
We will study separately each of the integrals $I_1$ and $I_2$ in the next lemmas:
\begin{Lemme}\label{Lemme41} For integral $I_1$ we have the estimate $I_1\leq C \mu\big(r+Ks_0\big)^{\omega-\gamma-1}$.
\end{Lemme}
\textit{\textbf{Proof}.} We begin by considering the space $\mathbb{R}^n$ as the union of a ball with dyadic coronas centered on $x(s_1)$, more precisely we set $\mathbb{R}^n=B_{f_1}\cup\bigcup_{k\geq 1}E_k$ where
\begin{eqnarray}\label{Decoupage}
B_{f_1}&=& \{x\in \mathbb{R}^n: |x-x(s_1)|\leq f_1\},\\[5mm]
E_k&=& \{x\in \mathbb{R}^n: f_12^{k-1}<|x-x(s_1)|\leq f_1 2^{k}\} \quad \mbox{for } k\geq 1.\nonumber
\end{eqnarray}
\begin{enumerate}
\item[(i)] \underline{Estimations over the ball $B_{f_1}$}. Applying Hölder's inequality on integral $I_1$ we obtain
\begin{eqnarray*}
I_{1,B_{f_1}}=\int_{B_{f_1}}|x-x(s_1)|^{\omega-1}|v-\overline{v}_{B_{f_1}}| |\psi_+(x,s_1)|dx & \leq &\underbrace{\||x-x(s_1)|^{\omega-1}\|_{L^p(B_{f_1})}}_{(1)} \\
& \times &\underbrace{\|v-\overline{v}_{B_{f_1}}\|_{L^z(B_{f_1})}}_{(2)}\underbrace{\|\psi_+(\cdot, s_1)\|_{L^q(B_{f_1})}}_{(3)}
\end{eqnarray*}
where $\frac{1}{p}+\frac{1}{z}+\frac{1}{q}=1$ and $p,z,q> 1$. We treat each of the previous terms separately:
\begin{enumerate}
\item[$\bullet$] Observe that for $1<p<n/(1-\omega)$ we have
\begin{equation*}
\||x-x(s_1)|^{\omega-1}\|_{L^p(B_{f_1})}\leq C f_1^{n/p+\omega-1}.
\end{equation*}
\item[$\bullet$] We have $v(\cdot, s_1)\in bmo(\mathbb{R}^n)$, thus $\|v-\overline{v}_{B_{f_1}}\|_{L^z(B_{f_1})}\leq C|B_{f_1}|^{1/z}\|v(\cdot,s_1)\|_{bmo}$. 
Since $\underset{s_1\in [s_0,T]}{\sup}\; \|v(\cdot,s_1)\|_{bmo}\leq \mu$ we write
\begin{equation*}
\|v-\overline{v}_{B_{f_1}}\|_{L^z(B_{f_1})}\leq C f_1^{n/z}\mu.
\end{equation*}
\item[$\bullet$] Finally, by the maximum principle for $L^q$ norms we have $\|\psi_+(\cdot, s_1)\|_{L^q(B_{f_1})}\leq  \|\psi(\cdot, s_0)\|_{L^q}$; hence we obtain
\begin{equation*}
\|\psi_+(\cdot, s_1)\|_{L^q(B_{f_1})}\leq  \|\psi(\cdot, s_0)\|_{L^1}^{1/q}\|\psi(\cdot, s_0)\|_{L^\infty}^{1-1/q}.
\end{equation*}
\end{enumerate}
We combine all these inequalities in order to obtain the following estimation for $I_{1,B_{f_1}}$:
\begin{equation*}
I_{1,B_{f_1}}\leq C\mu f_1^{n(1-1/q)+\omega-1}\|\psi(\cdot, s_0)\|_{L^1}^{1/q}\|\psi(\cdot, s_0)\|_{L^\infty}^{1-1/q}.
\end{equation*}
\item[(ii)] \underline{Estimations for the dyadic corona $E_k$}. Let us note $I_{1,E_k}$ the integral
$$I_{1,E_k}=\int_{E_k}|x-x(s_1)|^{\omega-1}|v-\overline{v}_{B_{f_1}}| |\psi_+(x,s_1)|dx.$$
Since over $E_k$ we have $|x-x(s_1)|^{\omega-1}\leq C 2^{k(\omega-1)}f_1^{\omega-1}$ we write
\begin{eqnarray*}
I_{1,E_k}&\leq & C2^{k(\omega-1)}f_1^{\omega-1}\left(\int_{E_k}|v-\overline{v}_{B(f_12^k)}| |\psi_+(x,s_1)|dx+\int_{E_k}|\overline{v}_{B_{f_1}}-\overline{v}_{B(f_12^k)}| |\psi_+(x,s_1)|dx\right)\\
&\leq& C2^{k(\omega-1)}f_1^{\omega-1}\left(\int_{B(f_12^k)}|v-\overline{v}_{B(f_12^k)}| |\psi_+(x,s_1)|dx\right.\\
& &\qquad \qquad \qquad \qquad\left.+\int_{B(f_12^k)}|\overline{v}_{B_{f_1}}-\overline{v}_{B(f_12^k)}| |\psi_+(x,s_1)|dx\right).
\end{eqnarray*}
where $B(f_12^k)= \{x\in \mathbb{R}^n: |x-x(s_1)|\leq f_12^k\}$.\\

Now, since $v(\cdot, s_1)\in bmo(\mathbb{R}^n)$, using the Lemma \ref{PropoBMO1} we have $|\overline{v}_{B_{f_1}}-\overline{v}_{B(f_12^k)}| \leq Ck\|v(\cdot,s_1)\|_{bmo}\leq Ck\mu$. We write 
\begin{eqnarray*}
I_{1,E_k}&\leq &C2^{k(\omega-1)}f_1^{\omega-1}\left(\int_{B(f_12^k)}|v-\overline{v}_{B(f_12^k)}| |\psi_+(x,s_1)|dx+Ck\mu\|\psi_+(\cdot,s_1)\|_{L^1}\right)\\[5mm]
&\leq & C2^{k(\omega-1)}f_1^{\omega-1}\left(\|\psi_+(\cdot,s_1)\|_{L^{a_0}} \|v-\overline{v}_{B(f_12^k)}\|_{L^{\frac{a_0}{a_0-1}}} +Ck\mu\;\|\psi_+(\cdot,s_0)\|_{L^1}\right)
\end{eqnarray*}
where we used Hölder's inequality with $1<a_0<\frac{n}{n+(\omega-1)}$ and maximum principle for the last term above. Using again the properties of $bmo$ spaces we have
$$I_{1,E_k}\leq  C2^{k(\omega-1)}f_1^{\omega-1}\left(\|\psi_+(\cdot,s_0)\|_{L^1}^{1/a_0}\|\psi_+(\cdot,s_0)\|_{L^\infty}^{1-1/a_0} |B(f_12^k)|^{1-1/a_0}\|v(\cdot,s_1)\|_{bmo} +Ck\mu\|\psi(\cdot,s_0)\|_{L^1}\right).$$
Since $\|v(\cdot,s_1)\|_{bmo}\leq \mu$ and since $1<a_0<\frac{n}{n+(\omega-1)}$, we have $n(1-1/a_0)+(\omega-1)<0$, so that, summing over each dyadic corona $E_k$, we obtain
\begin{equation*}
\sum_{k\geq 1}I_{1,E_k}\leq C \mu\left(f_1^{n(1-1/a_0)+\omega-1}\|\psi(\cdot,s_0)\|_{L^1}^{1/a_0}\|\psi(\cdot,s_0)\|_{L^\infty}^{1-1/a_0}+f_1^{\omega-1}\|\psi(\cdot,s_0)\|_{L^1}\right).
\end{equation*}
\end{enumerate}
We finally obtain the following inequalities:
\begin{eqnarray}
I_1&=& I_{1,B_{f_1}}+\sum_{k\geq 1}I_{1,E_k}\label{FinalEstimateI_1}\\
&\leq &C\mu \underbrace{f_1^{n(1-1/q)+\omega-1}\|\psi(\cdot, s_0)\|_{L^1}^{1/q}\|\psi(\cdot, s_0)\|_{L^\infty}^{1-1/q}}_{(a)}\nonumber\\
& &+ C \mu \left(\underbrace{f_1^{n(1-1/a_0)+\omega-1}\|\psi(\cdot,s_0)\|_{L^1}^{1/a_0}\|\psi(\cdot,s_0)\|_{L^\infty}^{1-1/a_0}}_{(b)} +\underbrace{f_1^{\omega-1}\|\psi(\cdot,s_0)\|_{L^1}}_{(c)}\right)\nonumber
\end{eqnarray}
Now we will prove that each of the terms $(a)$, $(b)$ and $(c)$ above is bounded by the quantity $\big(r+Ks_0\big)^{\omega-\gamma-1}$:
\begin{itemize}
\item for the first term (a) by the hypothesis on the initial data $\psi(\cdot,s_0)$ and the definition of $f_1$ given in (\ref{DefiFunctionF}) we have:
\begin{eqnarray*}
f_1^{n(1-1/q)+\omega-1}\|\psi(\cdot, s_0)\|_{L^1}^{1/q}\|\psi(\cdot, s_0)\|_{L^\infty}^{1-1/q}& \leq & \big(r+Ks_0\big)^{[n(1-1/q)+\omega-1]-\frac{\gamma}{q}-(n+\gamma)(1-1/q)}=\big(r+Ks_0\big)^{\omega-\gamma-1}.
\end{eqnarray*}
\item For the second term (b) we have, by the same arguments:
\begin{eqnarray*}
f_1^{n(1-1/a_0)+\omega-1}\|\psi(\cdot,s_0)\|_{L^1}^{1/a_0}\|\psi(\cdot,s_0)\|_{L^\infty}^{1-1/a_0} &\leq &
\big(r+Ks_0\big)^{[n(1-1/a_0)+\omega-1]-\frac{\gamma}{a_0}-(n+ \gamma)(1-1/a_0)}=\big(r+Ks_0\big)^{\omega-\gamma-1}.
\end{eqnarray*}
\item Finally, for the last term (c) we write
\begin{eqnarray*}
f_1^{\omega-1}\|\psi(\cdot,s_0)\|_{L^1}&\leq &f_1^{\omega-1}(r+Ks_0)^{-\gamma}= \big(r+Ks_0\big)^{\omega-\gamma-1}.
\end{eqnarray*}
\end{itemize}
Gathering these estimates on $(a), (b)$ and $(c)$, and getting back to (\ref{FinalEstimateI_1}) we finally obtain
$$I_1\leq C\mu \big(r+Ks_0\big)^{\omega-\gamma-1}.$$
The Lemma \ref{Lemme41} is proven.\hfill$\blacksquare$\\

\begin{Lemme}\label{Lemme42}
For integral $I_2$ in the inequality (\ref{Estrella}) we have the following estimate $I_2\leq C\big(r+Ks_0\big)^{\omega-\gamma-1}$.
\end{Lemme}
\textbf{\textit{Proof.}} As for the Lemma \ref{Lemme41}, we consider $\mathbb{R}^n$ as the union of a ball with dyadic coronas centered on $x(s_1)$ (cf. (\ref{Decoupage})). 
\begin{enumerate}
\item[(i)] \underline{Estimations over the ball $B_{f_1}$}. We will follow closely the computations of the Lemma \ref{Lemme2}. We write:
\begin{eqnarray*}
I_{2,B_{f_1}}&=&\int_{B_{f_1}}|\mathcal{L}(|x-x(s_1)|^{\omega})|\,|\psi_+(x,s_1)|dx \leq \|\psi_+(\cdot, s_1)\|_{L^{\infty}}\int_{B_{f_1}}|\mathcal{L}(|x-x(s_1)|^{\omega})|dx\\
&\leq & \|\psi_+(\cdot, s_0)\|_{L^\infty}\int_{\{|x|\leq f_1\}} \left|\mbox{v.p.}\int_{\mathbb{R}^n}[|x|^{\omega}-|x-y|^{\omega}]\pi(y)dy\right|dx.
\end{eqnarray*}
In the case \textbf{(c)} when $\alpha=\beta=1/2$ and $\delta<1/2$ we write:
\begin{eqnarray*}
I_{2,B_{f_1}} &\leq & \|\psi_+(\cdot, s_0)\|_{L^\infty}\left(\int_{\{|x|\leq f_1\}}\left|\mbox{v.p.}\int_{\{|y| \leq 1\}}\frac{|x|^{\omega}-|x-y|^{\omega}}{|y|^{n+1}}dy\right|dx+\int_{\{|x|\leq f_1\}} \int_{\mathbb{R}^n}\frac{||x|^{\omega}-|x-y|^{\omega}|}{|y|^{n+2\delta}}dydx\right)
\end{eqnarray*}
Following exactly the same arguments used in Lemma \ref{Lemme2} with the formulas (\ref{Estimate101})-(\ref{FirstHomogeneityEstimate}), \textit{i.e.} essentially by homogeneity, we have
$$I_{2,B_{f_1}}\leq  C\|\psi_+(\cdot, s_0)\|_{L^\infty} (f_1^{n+\omega-1}+f_1^{n+\omega-2\delta})$$
Since $0<2\delta <1$, recalling that by the definition of the function $f_1$ we have the estimate $0<f_1<1$, we obtain $f_1^{\omega-2\delta-\gamma}\leq f_1^{\omega-1-\gamma}$. 
The case \textbf{(d)} is straightforward since $\mathcal{L}=(-\Delta)^{1/2}$ and $(-\Delta)^{1/2}(|x|^{\omega})=|x|^{\omega-1}$.\\

Thus, in any case, we can write:
\begin{equation}\label{Bola2}
I_{2,B_{f_1}}\leq Cf_1^{n+\omega-1}\|\psi_+(\cdot, s_0)\|_{L^\infty}.
\end{equation}

\item[(ii)]  \underline{Estimations for the dyadic corona $E_k$}. Here we have
\begin{eqnarray*}
I_{2,E_k}=\int_{E_k}|\mathcal{L}(|x-x(s_1)|^{\omega})|\, |\psi_+(x,s_1)|dx & \leq &\|\psi_+(\cdot, s_0)\|_{L^1} \underset{f_12^{k-1}<|x|\leq f_12^k}{\sup}\left|\mbox{v.p.}\int_{\mathbb{R}^n}[|x|^{\omega}-|x-y|^{\omega}]\pi(y)dy\right|.
\end{eqnarray*}
In the case \textbf{(c)} we have:
\begin{eqnarray*}
I_{2,E_k}& \leq &\|\psi_+(\cdot, s_0)\|_{L^1}\underset{f_12^{k-1}<|x|\leq f_12^k}{\sup}\left( \left|\mbox{v.p.}\int_{\{|y|\leq 1\}}\frac{|x|^{\omega}-|x-y|^{\omega}}{|y|^{n+1}}dy\right|+\int_{\mathbb{R}^n}\frac{||x|^{\omega}-|x-y|^{\omega}|}{|y|^{n+2\delta}}dy\right).
\end{eqnarray*}
Again, by homogeneity and following the same lines of the Lemma \ref{Lemme2} above, we have
$$I_{2,E_k}\leq C\|\psi_+(\cdot, s_0)\|_{L^1}\left( (f_12^{k})^{\omega-1}\big(1+\ln(2^{k-1})\big)+(f_12^{k})^{\omega-2\delta}\right) $$
Since $0<\gamma<\omega<2\delta<1$ we have $\omega-1<0$ and $\omega-2\delta<0$ and thus, summing over $k\geq 1$, we obtain
\begin{equation*}
\sum_{k\geq 1}I_{2,E_k}\leq C\left( f_1^{\omega-1}+f_1^{\omega-2\delta}\right)\|\psi(\cdot, s_0)\|_{L^1}.
\end{equation*}
Repeating the same argument used before (\textit{i.e.} the fact that $0<f_1<1$), we finally get
\begin{equation}\label{Coronak2}
\sum_{k\geq 1}I_{2,E_k}\leq C f_1^{\omega-1}\|\psi(\cdot, s_0)\|_{L^1}.
\end{equation}
For the case \textbf{(d)}, we obtain the same inequality by homogeneity. 
\end{enumerate}
To finish the proof of the Lemma \ref{Lemme42} we combine (\ref{Bola2}) and (\ref{Coronak2}) and we obtain
$$I_2=I_{2,B_{f_1}}+\sum_{k\geq 1}I_{2,E_k}\leq C\left(\underbrace{f_1^{n+\omega-1}\|\psi_+(\cdot, s_0)\|_{L^\infty}}_{(d)}+\underbrace{f_1^{\omega-1}\|\psi(\cdot, s_0)\|_{L^1}}_{(e)}\right)$$
Now, we prove that the quantities $(d)$ and $(e)$ can be bounded by $\big(r+Ks_0\big)^{\omega-\gamma-1}$.
\begin{itemize}
\item For the term $(d)$ we write $f_1^{n+\omega-1}\|\psi(\cdot, s_0)\|_{L^\infty}\leq f_1^{n +\omega-1}(r+Ks_0)^{-(n+\gamma)}= \big(r+Ks_0\big)^{\omega-\gamma-1}$.
\item To treat the term $(e)$ it is enough to apply the same arguments used to prove the part $(c)$ above. 
\end{itemize}
Finally, we obtain
$$I_2=I_{2,B_{f_1}}+\sum_{k\geq 1}I_{2, E_k}\leq  C\big(r+Ks_0\big)^{\omega-\gamma-1}$$
The Lemma \ref{Lemme42} is proven.\hfill$\blacksquare$\\


Now we continue the proof of the Proposition \ref{Propo4}. Using the Lemmas \ref{Lemme41} and \ref{Lemme42} and getting back to the estimate (\ref{Estrella}) we have
\begin{equation}\label{FinalEstimate}
\left|\partial_{s_1} \int_{\mathbb{R}^n}\Omega_1(x)\psi_+(x,s_1)dx\right| \leq  C \left(\mu+1\right)\big(r+Ks_0\big)^{\omega-\gamma-1}
\end{equation}
This estimation is compatible with the estimate (\ref{Concentration2}) for $0\leq s_1\leq \epsilon r$ small enough. Indeed, we can write
$\phi=(r+K(s_0+s_1))^{\omega-\gamma}$ and we linearize this expression with respect to $s_1$:
$$\phi \thickapprox(r+s_0)^{\omega-\gamma}\left(1+K(\omega-\gamma)\frac{s_1}{(r+s_0)}\right)$$
Taking the derivative of $\phi$ with respect to $s_1$ we have $\phi' \thickapprox K(\omega-\gamma) \big(r+Ks_0\big)^{\omega-\gamma-1}$ and with the condition (\ref{SmallConstants}) on $K(\omega-\gamma)$ we obtain that (\ref{FinalEstimate}) is bounded by $\phi'$ and the Proposition \ref{Propo4} follows. \hfill$\blacksquare$\\

Now we write down the maximum principle for a small time $s_1$ but with a initial condition $\psi(\cdot,s_0)$, with $s_0>0$.
\begin{Proposition}[Height condition]\label{Propo5}
Under the  hypothesis of Theorem \ref{Generalisacion}, if $\psi(x,s_1)$ satisfies the concentration condition (\ref{Concentration2}), then we have the following height condition
\begin{equation*}
\|\psi(\cdot, s_1)\|_{L^\infty}  \leq \frac{1}{\left(r+K(s_0+s_1)\right)^{n+\gamma}}.
\end{equation*}
\end{Proposition}
\textit{\textbf{Proof.}} The proof follows essentially the same lines of the Proposition \ref{Propo2}. Indeed, since we have assumed that the concentration condition (\ref{Concentration2}) is bounded by $(r+K(s_0+s_1))^{\omega-\gamma}$, we obtain in the same manner and with the same constants:
$$\frac{d}{ds_1}\|\psi(\cdot,s_1)\|_{L^\infty}\leq -C(r+K(s_0+s_1))^{-\frac{(\omega-\gamma)}{n+\omega}}\|\psi(\cdot,s_1)\|_{L^\infty}^{1+\frac{1}{n+\omega}}.$$
To conclude, it is enough to solve the previous differential inequality with initial data $\|\psi(\cdot,0)\|_{L^\infty}\leq (r+Ks_0)^{-(n+\gamma)}$ to obtain that $\|\psi(\cdot,s_1)\|_{L^\infty}\leq (r+K(s_0+s_1))^{-(n+\gamma)}$. \hfill$\blacksquare$\\

The crucial part of the proof of Theorem \ref{Generalisacion} is given by the next proposition which gives us a control on the $L^1$-norm for a time $s_0+s_1$.
\begin{Proposition}[Second $L^1$-norm estimate]\label{Propo6}
Under the hypothesis of Theorem (\ref{Generalisacion}) we have
$$\|\psi(\cdot,s_1)\|_{L^1}\leq \frac{v_n}{\big(r+K(s_0+s_1)\big)^{\gamma}}$$
\end{Proposition}
\textit{\textbf{Proof.}} This is a direct consequence of the concentration condition and of the previous height condition.
\hfill$\blacksquare$
\subsection{The iteration}\label{SecIterationMol}

In sections \ref{SecEvolMol1} and \ref{SecEvolMol2} we studied respectively the evolution of small molecules from time $0$ to a small time $s_0$ and from this time $s_0$ to a larger time $s_0+s_1$ and we obtained a good $L^1$ control for such molecules. It is now possible to reapply the previous Theorem \ref{Generalisacion} in order to obtain a larger time control of the $L^1$ norm. The calculus of the $N$-th iteration will be essentially the same.
\begin{Theoreme}\label{Generalisacionfin} Set $\gamma$ and $\omega$ two real numbers such that $0<\gamma<\omega<2\delta<1$ in the case \textbf{(c)} or $0<\gamma<\omega<1$ in the case \textbf{(d)}. Let $0< s_N\leq T$ and let $\psi(x,s_N)$ be a solution of the problem
\begin{equation}\label{Evolutionfin}
\left\lbrace
\begin{array}{rl}
\partial_{s_N} \psi(x,s_N)=& -\nabla\cdot(v\, \psi)(x,s_N)-\mathcal{L}\psi(x,s_N)\\[5mm]
\psi(x,0)=& \psi(x,s_{N-1})  \qquad \qquad \mbox{with } s_{N-1}>0\\[5mm]
div(v)=&0 \quad \mbox{and }\; v\in L^{\infty}([0,T];bmo(\mathbb{R}^n))\quad \mbox{with } \underset{s_N\in [s_{N-1},T]}{\sup}\; \|v(\cdot,s_N)\|_{bmo}\leq \mu
\end{array}
\right.
\end{equation}
If $\psi(x,s_{N-1})$ satisfies the three following conditions
\begin{eqnarray*}
\int_{\mathbb{R}^n}|\psi(x,s_{N-1})||x-x(s_{N-1})|^\omega dx &\leq &(r+K(s_0+\cdots+s_{N-1}))^{\omega-\gamma}\\
\|\psi(\cdot, s_{N-1})\|_{L^\infty}\leq  \frac{1}{\left(r+ K(s_0+\cdots+s_{N-1})\right)^{n+\gamma }}&; &\quad \|\psi(\cdot, s_{N-1})\|_{L^1} \leq  \frac{v_n}{\big(r+K(s_0+\cdots+s_{N-1})\big)^{\gamma }}
\end{eqnarray*}
where $K=K(\mu)$ is given by (\ref{SmallConstants}) and $s_N$ is such that $(r+K(s_0+\cdots+s_N))<1$. Then for all $0< s_N\leq\epsilon r$ small, we have the following estimates
\begin{eqnarray}
\int_{\mathbb{R}^n}|\psi(x,s_N)||x-x(s_N)|^\omega dx &\leq &(r+K(s_0+\cdots+s_N))^{\omega-\gamma}  \label{Concentration2fin}\\
\|\psi(\cdot,s_N)\|_{L^\infty}&\leq & \frac{1}{\left(r+K(s_0+\cdots+s_N)\right)^{n+\gamma }}\nonumber\\
\|\psi(\cdot,s_N)\|_{L^1} &\leq & \frac{v_n}{\big(r+K(s_0+\cdots+s_N)\big)^{\gamma }}\nonumber
\end{eqnarray}
\end{Theoreme}

\begin{Remarque}
\begin{itemize}
\item[]
\item[1)] Again, since $s_N$ is small and $(r+K(s_0+\cdots+s_{N-1}))<1$, we can without loss of generality assume that $(r+K(s_0+\cdots +s_N))<1$: otherwise, by the maximum principle there is nothing to prove.
\item[2)] The new molecule's center $x(s_N)$ used in formula (\ref{Concentration2fin}) is fixed by 
\begin{equation}\label{Defpointx_sfin}
\left\lbrace
\begin{array}{rl}
x'(s_N)=& \overline{v}_{B_{f_N}}=\frac{1}{|B_{f_N}|}\displaystyle{\int_{B_{f_N}}}v(y,s_N)dy\\[5mm]
x(0)=& x(s_{N-1}).
\end{array}
\right.
\end{equation}  
And here we noted $B_{f_N}=B(x(s_N),f_N)$ with $f_N$ a real valued function given by 
\begin{equation}\label{DefiFunctionFfin}
f_N=(r+K(s_0+\cdots+s_{N-1})).
\end{equation}
Note that by remark $1)$ above we have $0<f_N<1$.
\end{itemize}
\end{Remarque}

The proof of Theorem \ref{Generalisacionfin} will follow exactly the same steps given in the proof of Theorem \ref{Generalisacion}: we start with the concentration condition studied in Proposition \ref{Propo4fin} and we continue with the Height condition in Proposition \ref{Propo5fin}, finally, the $L^1$ bound will be an easy consequence of these two estimates.
\begin{Proposition}[Concentration condition]\label{Propo4fin}
Under the  hypothesis of Theorem \ref{Generalisacionfin}, if $\psi(\cdot,s_{N-1})$ is an initial data then the solution $\psi(x,s_N)$ of (\ref{Evolutionfin}) satisfies
\begin{equation*}
\int_{\mathbb{R}^n} |\psi(x,s_N)||x-x(s_N)|^{\omega}dx \leq (r+K(s_0+\cdots+s_N))^{\omega-\gamma}
\end{equation*}
for $x(s_N)\in \mathbb{R}^n$ fixed by formula (\ref{Defpointx_sfin}), with $0\leq s_N\leq \epsilon r$. 
\end{Proposition}
\textit{\textbf{Proof}.} Follow the same lines given in the proof of Proposition \ref{Propo4}. Write $\Omega_N(x)=|x-x(s_N)|^{\omega}$ and $\psi(x)=\psi_+(x)-\psi_-(x)$, by linearity and using the positivity principle we have  $|\psi(x,s_N)|=|\psi_+(x,s_N)-\psi_-(x,s_N)|\leq \psi_+(x,s_N)+\psi_-(x,s_N)$ and we may consider the formula:
\begin{eqnarray*}
I&=&\left|\partial_{s_N} \int_{\mathbb{R}^n}\Omega_N(x)\psi_+(x,s_N)dx\right|=\left|\int_{\mathbb{R}^n}-\nabla\Omega_N(x)\cdot x'(s_N)\psi_+(x,s_N)+\Omega_N(x)\left[-\nabla\cdot(v\, \psi_+(x,s_N))-\mathcal{L}\psi_+(x,s_N)\right]dx\right|
\end{eqnarray*}
Using the definition of $x'(s_N)$ given in (\ref{Defpointx_sfin}) and replacing $\Omega_N(x)$ by $|x-x(s_N)|^{\omega}$ in the first integral we obtain
\begin{equation}\label{Estrellafin}
I\leq c \underbrace{\int_{\mathbb{R}^n}|x-x(s_N)|^{\omega-1}|v-\overline{v}_{B_f}| |\psi_+(x,s_N)|dx}_{I_1} + c\underbrace{\int_{\mathbb{R}^n}|\mathcal{L}\Omega_N(x)||\psi_+(x,s_N)|dx}_{I_2}.
\end{equation}
We will study each of the integrals $I_1$ and $I_2$ in the next lemmas:
\begin{Lemme}\label{Lemme41fin} For integral $I_1$ we have $I_1\leq C \mu\big(r+K(s_0+\cdots+s_{N-1})\big)^{\omega-\gamma-1}$.
\end{Lemme}
\textit{\textbf{Proof of the lemma}.} It is enough to repeat the same steps of the previous Lemma \ref{Lemme41}, just consider $\mathbb{R}^n=B_{f_N}\cup\bigcup_{k\geq 1}E_k$ where
\begin{eqnarray}\label{Decoupagefin}
B_{f_N}= \{x\in \mathbb{R}^n: |x-x(s_N)|\leq f_N\}, \qquad E_k= \{x\in \mathbb{R}^n: f_N2^{k-1}<|x-x(s_N)|\leq f_N 2^{k}\} \quad \mbox{for } k\geq 1.
\end{eqnarray}
In order to obtain the desired inequality, use exactly the same arguments, the maximum principle and the hypothesis of Theorem  \ref{Generalisacionfin}.\hfill$\blacksquare$
\begin{Lemme}\label{Lemme42fin}
For integral $I_2$ in inequality (\ref{Estrellafin}) we have the following estimate
\begin{equation*}
I_2=\int_{\mathbb{R}^n}|\mathcal{L}\Omega_N(x)||\psi_+(x,s_N)|dx\leq C\big(r+K(s_0+\cdots+s_{N-1})\big)^{\omega-\gamma-1}.
\end{equation*}
\end{Lemme}
\textbf{\textit{Proof of the lemma.}} As for Lemma \ref{Lemme41fin}, we consider $\mathbb{R}^n$ as the union of a ball with dyadic coronas centered on $x(s_N)$ (cf. (\ref{Decoupagefin})). It is then enough to repeat the corresponding estimates of the $s_1$-case given in Lemma \ref{Lemme42}. \hfill$\blacksquare$\\

Now we continue the proof of the Proposition \ref{Propo4fin}. Using the Lemmas \ref{Lemme41fin} and \ref{Lemme42fin} and getting back to the estimate (\ref{Estrellafin}) we have
\begin{equation}\label{FinalEstimatefin}
\left|\partial_{s_N} \int_{\mathbb{R}^n}\Omega_N(x)\psi_+(x,s_N)dx\right| \leq  C \left(\mu+1\right)\big(r+K(s_0+\cdots+s_{N-1})\big)^{\omega-\gamma-1}
\end{equation}
This estimation is compatible with the estimate (\ref{Concentration2fin}) for $0\leq s_N\leq \epsilon r$ small enough. Indeed, we can write
$\phi=(r+K(s_0+\cdots+s_N))^{\omega-\gamma}$ and we linearize this expression with respect to $s_N$:
$$\phi \thickapprox(r+K(s_0+\cdots+s_{N-1}))^{\omega-\gamma}\left(1+K(\omega-\gamma)\frac{s_N}{(r+K(s_0+\cdots+s_{N-1}))}\right)$$
Taking the derivative of $\phi$ with respect to $s_N$ we have $\phi' \thickapprox K(\omega-\gamma) \big(r+K(s_0+\cdots+s_{N-1})\big)^{\omega-\gamma-1}$ and with the condition (\ref{SmallConstants}) on $K(\omega-\gamma)$ we obtain that (\ref{FinalEstimatefin}) is bounded by $\phi'$ and the Proposition \ref{Propo4fin} follows. \hfill$\blacksquare$
\begin{Proposition}[Height condition]\label{Propo5fin}
Under the  hypothesis of Theorem \ref{Generalisacionfin}, if $\psi(x,s_N)$ satisfies concentration condition (\ref{Concentration2fin}), then we have the next height condition
\begin{equation*}
\|\psi(\cdot, s_N)\|_{L^\infty}  \leq \frac{1}{\left(r+K(s_0++\cdots+s_N)\right)^{n+\gamma}}.
\end{equation*}
\end{Proposition}
\textit{\textbf{Proof.}} The proof follows essentially the same lines of the Proposition \ref{Propo2}. Indeed, since we have that concentration condition (\ref{Concentration2fin}) is bounded by $(r+K(s_0+\cdots+s_N))^{\omega-\gamma}$, we obtain in the same manner and with the same constants:
$$\frac{d}{ds_N}\|\psi(\cdot,s_N)\|_{L^\infty}\leq -C(r+K(s_0+\cdots+s_N))^{-\frac{(\omega-\gamma)}{n+\omega}}\|\psi(\cdot,s_N)\|_{L^\infty}^{1+\frac{1}{n+\omega}}.$$
Solving this differencial inequality we obtain $\|\psi(\cdot,s_N)\|_{L^\infty}\leq (r+K(s_0+\cdots+s_N))^{-(n+\gamma)}$.\hfill$\blacksquare$
\begin{Proposition}[$L^1$-norm estimate]\label{Propo6fin}
Under the hypothesis of Theorem \ref{Generalisacionfin} we have
\begin{equation}\label{FinalEstimateFin}
\|\psi(\cdot,s_N)\|_{L^1}\leq \frac{v_n}{\big(r+K(s_0+\cdots+s_N)\big)^{\gamma}}
\end{equation}
\end{Proposition}
\textit{\textbf{Proof.}} This is a direct consequence of the concentration condition and of the previous height condition.
\hfill$\blacksquare$\\

\textit{\textbf{End of the proof of Theorem \ref{TheoL1control}.}} We have proved with the Theorem \ref{SmallGeneralisacion} that is possible to control the $L^1$ behavior of the molecules $\psi$ from $0$ to a small time $s_0$, from time $s_0$ to time $s_1$ with Theorem \ref{Generalisacion}, and by iteration from time $s_{N-1}$ to time $s_N$ with Theorem \ref{Generalisacionfin}. We recall that we have $s_i\sim \epsilon r$ for all $0\leq i\leq N$, so the bound obtained in (\ref{FinalEstimateFin}) depends mainly on the size of the molecule $r$ and the number of iterations $N$. \\

We observe now that the smallness of $r$ and of the times $s_0,...,s_N$ can be compensated by the number of iterations $N$ in the following sense: fix a small $0<r<1$ and iterate as explained before. Since each small time $s_0,...,s_N $ is of order $\epsilon r$, we have $s_0+\cdots+s_{N}\sim N \epsilon r$. Thus, we will stop the iterations as soon as $ N r\geq T_0$. 

Of course, the number of iterations $N=N(r)$ will depend on the smallness of the molecule's size $r$, and more specifically it is enough to set $N(r)\sim \frac{T_0}{r}$ in order to obtain this lower bound for $Nr$.

Proceeding this way we will obtain $\|\psi(\cdot,s_N)\|_{L^1}\leq C T_0^{-\gamma}<+\infty$, for all molecules of size $r$. Note in particular that, once this estimate is available, for bigger times it is enough to apply the maximum principle.\\

Finally, and for all $r>0$, we obtain after a time $T_0$ a $L^1$ control for small molecules and we finish the proof of the Theorem \ref{TheoL1control}. \hfill$\blacksquare$

\quad\\[5mm]

\begin{flushright}
\begin{minipage}[r]{80mm}
Diego \textsc{Chamorro}\\[3mm]
Laboratoire d'Analyse et de Probabilités\\ 
Université d'Evry Val d'Essonne\\[2mm]
23 Boulevard de France\\
91037 Evry Cedex\\[2mm]
diego.chamorro@univ-evry.fr
\end{minipage}
\end{flushright}

\end{document}